\documentclass[journal]{IEEEtran}
\usepackage[utf8]{inputenc}

\title{Planning Low-Carbon Campus Energy Hubs}
\author{Daniel~J.~Olsen,~\IEEEmembership{Student~Member,~IEEE},
        Ning~Zhang,~\IEEEmembership{Senior~Member,~IEEE},
        Chongqing~Kang,~\IEEEmembership{Fellow,~IEEE},
        Miguel~A.~Ortega-Vazquez,~\IEEEmembership{Senior~Member,~IEEE},
        Daniel~S.~Kirschen,~\IEEEmembership{Fellow,~IEEE}
}% <-this % stops a space}

\usepackage{amsmath}
\usepackage[table]{xcolor}  % http://ctan.org/pkg/xcolor
\usepackage{graphicx}
\usepackage{gensymb}        % \degree
\usepackage{cite}           %condenses [5], [6], [7] to [5-7]

\usepackage{tikz}
\usepackage{textcomp}
\usepackage{hyperref}
\usepackage{lipsum}

\newcommand{\copyrighttext}{%
    \footnotesize \textcopyright 2018 IEEE. Personal use of this material is permitted. Permission from IEEE must be obtained for all other uses, in any current or future media, including reprinting/republishing this material for advertising or promotional purposes, creating new collective works, for resale or redistribution to servers or lists, or reuse of any copyrighted component of this work in other works.	DOI: \href{https://ieeexplore.ieee.org/document/8523813}{10.1109/TPWRS.2018.2879792}, IEEE Transactions on Power Systems.}

\newcommand{\copyrightnotice}{%
    \begin{tikzpicture}[remember picture,overlay]
    \node[anchor=south,yshift=8pt] at (current page.south) {\fbox{\parbox{\dimexpr\textwidth-\fboxsep-\fboxrule\relax}{\copyrighttext}}};
    \end{tikzpicture}%
}
\newenvironment{ldescription}[1]
  {\begin{list}{}%
   {\renewcommand\makelabel[1]{##1\hfill}%
   \settowidth\labelwidth{\makelabel{#1}}%
   \setlength\leftmargin{\labelwidth}
   \addtolength\leftmargin{\labelsep}}}
  {\end{list}}

\usepackage{stackengine}
\def\yenrule{\rule{1.3ex}{.1ex}}
\def\textyen{\renewcommand\stacktype{L}\stackon[.4ex]{\stackon[.65ex]{Y}{\yenrule}}{\yenrule}}

\begin{document}

\maketitle
\copyrightnotice

\begin{abstract}
    Multi-energy systems can provide a constant level of service to end-use energy demands, while deriving delivered energy from a variety of primary/secondary energy sources. This fuel-switching capability can be used to reduce operating expenses, reduce environmental impacts, improve flexibility to accommodate renewable energy, and improve reliability.
    
    This paper presents four frameworks for incentivizing energy hub equipment investments for low-carbon operation targets. These frameworks vary in the measures taken to achieve low-carbon operation (explicit constraint vs. carbon pricing) and in the relationship between the hub builder and operator (cooperative vs. uncoordinated). The underlying energy hub model upon which these frameworks are built is an enhanced greenfield model, introducing `energy buses' to reduce dimensionality.
    
    A case study is conducted for a campus being designed in Beijing, and results from each framework are compared to illustrate their relative costs. When the operator cannot be trusted to cooperate in controlling emissions, the system must be `overbuilt' with more expensive equipment to ensure emissions target are met. A taxation-based approach increases overall costs at moderate emissions targets, but this effect decreases at aggressive targets. This paper also compares the cost of less efficient institutional frameworks with the most efficient approach, i.e. cooperation between builder and operator with constraints on emissions.
\end{abstract}

%\begin{IEEEkeywords}
%Energy hub, multiple-energy system, climate policy, carbon taxation, strategic planning, power system economics.
%\end{IEEEkeywords}

\vspace{-10pt}   %condense
\section*{Nomenclature}

\let\thefootnote\relax\footnote{This work was supported in part by the National Science Foundation of China (No. 51620105007).}
\let\thefootnote\relax\footnote{D. J. Olsen and D. S. Kirschen are with the University of Washington Department of Electrical Engineering in Seattle, WA 98195 USA (e-mail: \{djolsen, kirschen\}@uw.edu).}
\let\thefootnote\relax\footnote{N. Zhang, C. Kang are with the State Key Lab of Power Systems, Department of Electrical Engineering, Tsinghua University, Beijing 100084, China (e-mail: ningzhang@tsinghua.edu.cn).}
\let\thefootnote\relax\footnote{M. A. Ortega-Vazquez is with Grid Operations and Planning, Electric Power
Research Institute, Palo Alto, CA 94304 USA (e-mail: mortegavazquez@epri.com).}

\vspace{-20pt}   %condense
\subsection*{Abbreviations and Symbols}
\begin{ldescription}{$xxxxx$}
    \item [$\text{tCO}_{2}\text{e}$] Metric tons of greenhouse gases (GHG) converted to $\text{CO}_{2}$ equivalent.
    \item [\textyen] Chinese yuan (RMB).
\end{ldescription}

\subsection*{Sets and indices}

\begin{ldescription}{$xxxxx$}
    \item [$G$] Set of equipment (power conversion and storage devices), indexed by $g$.
    \item [$G^{\text{C}}$] Set of power conversion devices ($G^{\text{C}} \subset G$).
    \item [$G^{\text{S}}$] Set of storage devices ($G^{\text{S}} \subset G$).
    \item [$L$] Set of branch flows, indexed by $l$.
    \item [$M$] Set of energy types, indexed by $m$.
    \item [$N^{\text{A}}$] Set of input power capacity discretization binaries, indexed by $n^a$.
    \item [$N^{\text{B}}$] Set of storage energy discretization binaries, indexed by $n^b$.
    \item [$N^{\text{C}}$] Set of storage power discretization binaries, indexed by $n^c$.
    \item [$N^{\text{D}}$] Set of converter count discretization binaries, indexed by $n^d$.
    \item [$P$] Set of equipment ports, indexed by $p$.
    \item [$P^{\text{out}}$] Set of equipment output ports ($P^{\text{out}} \subset P$).
    \item [$P_{g}$] Set of ports of equipment $g$ ($P_{g} \subset P$).
    \item [$S$] Set of representative days, indexed by $s$.
    \item [$T$] Set of time periods, indexed by $t$.
    \item [$Y$] Set of years, indexed by $y$.
\end{ldescription}

\subsection*{Topology Matrices}

\begin{ldescription}{$xxxxx$}
    \item [$\boldsymbol{A}$] Network topology matrix, dimension ($P \times L$). When subscripted, $\boldsymbol{A}_g$ refers to the rows of the matrix corresponding to the ports of equipment $g$.% ($P_g \times L$).
    \item [$\boldsymbol{H}$] Converter efficiency matrix, dimension ($P^{\text{out}} \times P$). When subscripted, $\boldsymbol{H}_g$ refers to the rows of the matrix corresponding to the output ports of equipment $g$.% ($P_g \times L$).
    \item [$\boldsymbol{J}$] `Limiting' port matrix, dimension ($G \times P$).
    \item [$\boldsymbol{K}$] Branch directionality vector, dimension ($L$).
    \item [$\boldsymbol{U}$] Input port matrix, dimension ($M \times L$).
    \item [$\boldsymbol{W}$] Output port matrix, dimension ($M \times L$).
    \item [$\boldsymbol{Z}$] Network efficiency matrix, dimension ($P^{\text{out}} \times L$).
\end{ldescription}
Topology matrix values are defined in Appendix A.

%\vspace{-4pt}
\subsection*{Fixed Parameters}

\begin{ldescription}{$xxxxx$}
    \item [$C_{g}^{\text{unit}}$] Cost of one piece of equipment of energy conversion device $g$ (\textyen).
    \item [$C_{g}^{\text{power}}$] Per-unit cost of power for storage device $g$ (\textyen /MW).
    \item [$C_{g}^{\text{energy}}$] Per-unit cost of energy for storage device $g$ (\textyen /MWh).
    \item [$C_{m}^{\text{cap}}$] Cost of input capacity for energy $m$ (\textyen /MW).
    \item [$E^{\text{max}}$] Annual GHG emissions limit ($\text{tCO}_{2}\text{e}$).
    \item [$i$] Discount rate for calculating net-present value.
    \item [$\Delta t$] Time interval length (hours).
    \item [$\pi_{s}$] Probability of representative day $s$.
\end{ldescription}

\subsection*{Time-Varying Energy Parameters}

\noindent Subscripted ${}_{m,s,t,y}$ for energy $m$ at time $t$ on day $s$ in year $y$.
\begin{ldescription}{$xxxxx$}
    \item [$\boldsymbol{B}$] Power availability of input energy relative to its input capacity. For electricity generated by renewable sources, $0 \le B_{m,s,t,y} \le 1$. For seasonally available energy flows (\textit{e.g.} district heating), $B_{m,s,t,y} \in \left \{0,1 \right \}$. For all other energy flows, $B_{m,s,t,y} = 1$.
    \item [$\boldsymbol{e}$] Marginal emissions rate ($\text{tCO}_{2}\text{e}$/MWh).
    \item [$\boldsymbol{f}$] Input energy price (\textyen /MWh).
    \item [$\boldsymbol{h}$] Grid energy feed-in price (\textyen /MWh).
    \item [$\boldsymbol{L}$] End-use power demand (MW).
\end{ldescription}

\subsection*{Carbon Pricing Variables}

\begin{ldescription}{$xxxxx$}
    \item [$P^{\text{CO}_{2}}$] Tax rate for GHG emissions (\textyen /$\text{tCO}_{2}\text{e}$).
    \item [$SCoC$] Social cost of carbon rate (\textyen /$\text{tCO}_{2}\text{e}$).
\end{ldescription}

\subsection*{Investment Variables}

\begin{ldescription}{$xxxxx$}
    \item [$C^{\text{invest}}$] Total investment cost (\textyen).
    \item [$D_{g}^{\text{max}}$] Rated power for equipment $g$, for one piece of equipment for conversion devices, or total capacity for storage devices (MW).
    \item [$I_{g}$] Number of pieces of equipment purchased for energy conversion device $g$.
    \item [$P_{m}^{\text{max}}$] Purchased input power capacity for energy $m$ (\textyen).
    \item [$Q_{g}^{\text{max}}$] Purchased energy capacity for storage device $g$ (MWh).
\end{ldescription}

\subsection*{Operational Variables}

\noindent Hub operations during each representative day $s$ in each year $y$ are independent, so these indices are omitted when possible for brevity.
\begin{ldescription}{$xxxxx$}
    \item [$C^{\text{operate}}$] Hub operating cost (\textyen).
    \item [$E^{\text{operate}}$] GHG emissions from hub operation  ($\text{tCO}_{2}\text{e}$).
    \item [$P_{m,t}$] Power flow into the energy hub from the grid for energy $m$ at time $t$ (MW).
    \item [$Q_{g,t}$] State of charge of storage device $g$ at time $t$ (MWh).
    \item [$r_{m,t}$] Power flow out of the energy hub to the grid for energy $m$ at time $t$ (MW).
    \item [$T^{\text{operate}}$] Tax bill for GHG emissions (\textyen).
    \item [$V_{l,t}$] Power flow within the energy hub for branch $l$ at time $t$ (MW).
\end{ldescription}

\subsection*{Dual Variables}

\begin{ldescription}{$xxxxxxxx$}
    \item [$\alpha_{p,t}$] Dual variable for network power balance constraints at port $p$ during time $t$.
    \item [$\beta_{g,t}$] Dual variable for conservation of energy constraints for storage device $g$ at time $t$.
    \item [$\underline{\gamma}_{g,t}, \overline{\gamma}_{g,t}$] Dual variables for \{lower, upper\} state of charge constraint for storage device $g$ at time $t$.
    \item [$\zeta_{g,t}$] Dual variable for maximum power constraints for converter $g$ at time $t$.
    \item [$\underline{\kappa}_{g,t}, \overline{\kappa}_{g,t}$] Dual variables for \{lower, upper\} power constraint for storage device $g$ at time $t$.
    \item [$\rho_{m,t}$] Dual variable for input capacity constraints for energy flow $m$ at time $t$.
    \item [$\mu_{m,t}$] Dual variable for hub outflow constraints for energy flow $m$ at time $t$.
    \item [$\underline{\phi}_{m,t}, \overline{\phi}_{m,t}$] Dual variables for \{lower, upper\} grid sales constraint for energy flow $m$ at time $t$.
    \item [$\sigma_{l,t}$] Dual variable for branch flow non-negativity constraint for branch $l$ at time $t$.
\end{ldescription}

\section{Introduction}

\IEEEPARstart{I}{nternational} efforts to cooperatively address anthropogenic climate change \cite{ipcc_syr_2014} have resulted in the Paris Climate Accords (PCA), in which each signatory country has made pledges to reduce their emissions of greenhouse gases (GHG) \cite{carbonbrief2017}. However, the current pledges are projected to result in a temperature rise of 3.2\degree C over pre-industrial levels by the year 2100, above the 2\degree C goal set in the PCA. Worse still, the current policies are projected to result in a temperature rise of 3.4\degree C \cite{climate_action_tracker}. In order to limit the worst effects of climate change, more aggressive regulatory policies toward reducing GHG emissions are necessary. Since building-induced GHG emissions are a significant source of global emissions, they are a prime candidate for emissions reduction policies.

Buildings and their occupants require several forms of energy (e.g. electrical, thermal, kinetic), each of which can often be supplied via several means. For example, thermal energy can be delivered via district heating or via conversion from electricity or natural gas. By considering all energy requirements and equipment of a multiple-energy system simultaneously, overall operation can be improved; this modeling approach is often referred to as using \textit{energy hubs} \cite{geidl2007}, which create algebraic representations of energy flows into, out of, and within a system for the purpose of optimization. Operational improvements can be measured in terms of cost, reliability, environmental impact, or other metrics. For example, an energy hub may have electricity and natural gas as inputs, require electricity and heat as output, and have a heat pump and a combined-heat-and-power (CHP) unit as energy conversion equipment. An optimization for cost may result in sourcing heat primarily from natural gas via the CHP unit, while an optimization for local air-quality may result sourcing heat primarily form electricity via the heat pump. Optimization of energy hub scheduling considering the Pareto frontier of cost/emissions tradeoffs is presented in \cite{schwaegerl2011}.

Because energy flows in an energy hub are constrained by planning and construction decisions as well as the selection of equipment type and capacity, these decisions have a long-lasting impact. The design of an entirely new multiple-energy system is sometimes known as \textit{greenfield} design. Energy hub models for the optimization of planning and operation of multiple-energy systems have been developed at the building \cite{bozchalui2012, evins2015, yi2018, picard2018}, district \cite{weber2011, morvaj2016, wang2017apen}, and regional \cite{chaudry2014} scales. An energy hub framework which incorporates building equipment and plug-in hybrid electric vehicle charging to respond to frequency control is presented in \cite{galus2011}. A review of energy hub models and other related modeling frameworks (\textit{e.g.} microgrids) can be found in \cite{mancarella2014}.

Several papers have investigated simultaneous optimization of planning and operation in energy hubs while considering carbon emissions. Considering carbon emissions as an additional objective to be minimized creates a Pareto frontier of solutions which trade off emissions reductions for cost increases, and vice versa. The full Pareto frontier of cost and emissions can be sampled using the $\epsilon\text{-constraint}$ method \cite{yi2018, weber2011, morvaj2016}, or if the set of potential technologies is small, the results for all combinations can be explicitly calculated \cite{picard2018}. Evins \cite{evins2015} presents a multi-level model where building and energy hub variables are optimized in the upper-level while operation variables (including binary variables for fuel cell status) are optimized in the lower level. Solutions are found using a non-dominated sorting genetic algorithm in order to trace out the Pareto frontier.

Another approach to sampling the Pareto frontier is a linear weighting factor, which can find points on the convex hull of the frontier. In this context, the weighting factor is a price for carbon emissions. Regulatory imposition of a price on externalities improves overall social welfare \cite{pigou}, and is one approach to avoiding a `tragedy of the commons' outcome where individually rational decisions result in a socially-suboptimal solution compared to a cooperative approach \cite{hardin1968}. Carbon pricing is increasingly common \cite{worldbank2017}, and can be implemented via a real price (a tax or an emissions-trading scheme) or a mandate to consider the \textit{social cost of carbon} (SCoC) in planning decisions \cite{tsd2010, tsd2013_2016, nordhaus2017}.

The motivation for this paper is to bridge perceived gaps in the existing literature: a), a mixed-inter linear program (MILP) model for planning and incentivizing low-carbon energy hubs, considering an independent operator that may not share the low-carbon goals of planners, and b) incentivization of low-carbon goals via an optimized price for carbon. Picard and Helsen \cite{picard2018} evaluate only a limited number of possible equipment combinations in order to be able to evaluate all of their costs and emissions. Several authors \cite{yi2018, weber2011, morvaj2016} build MILP models, but without an independent operator. Therefore the derived solutions \textit{could} satisfy the emissions constraints, but may not in practice. Evins \cite{evins2015} incorporates an independent, emissions-indifferent operator, but their approach utilizes an external building simulation (EnergyPlus) and candidate solutions are found via a genetic algorithm, so the quality of the best currently-found solution relative to the true optimum is unknown.

This paper extends the analysis of low-carbon energy hub design to include two strategic considerations. The first makes investment decisions while accounting for a hub operator that may ignore emissions-reduction goals, and the second decides carbon prices to induce lower-emission investment and operation decisions. The underlying greenfield energy hub model is also enhanced.
Specifically, this paper makes the following contributions:

\begin{itemize}
    \item The formulation of four different frameworks for optimizing low-carbon energy hub investment and operation, to account for differences in policy and market structures: Single Builder-Operator, Bi-level Regulator/Builder-Operator, Bi-level Builder/Operator, and Tri-level Regulator/Builder/Operator.
    \item An expanded investment optimization problem, including on-site renewable generation, grid capacity costs, and storage for each energy type.
    \item A new formulation of the greenfield energy hub operation model using the concept of \textit{energy buses}. This formulation simplifies the investment optimization problem without loss of accuracy.
\end{itemize}

\section{Low-Carbon Design Frameworks} \label{design_frameworks}

A strategy for controlling operational carbon emissions when planning equipment investments for a greenfield energy hub depends on the answers to two fundamental questions:
\begin{enumerate}
    \item Does the hub operator share the hub builder's goal to reduce carbon emissions?
    \item Are carbon emissions controlled via an explicit constraint or a carbon price set to meet a given target?
\end{enumerate}
The answers to these questions determine the four possible frameworks defined below and illustrated in Fig. \ref{frameworks}.

\begin{itemize}
    \item \textbf{Framework 1: Single Builder-Operator, Emission-Constrained}: A single entity designs and operates the hub at or below a given emissions target.
    \item \textbf{Framework 2: Bi-Level Regulator/Builder-Operator, Carbon Tax}: A regulatory agency sets a carbon tax rate, such that a builder-operator's minimum-cost investment and operation solution results in emissions at or below a given emissions target.
    \item \textbf{Framework 3: Bi-Level Builder/Operator, Emission Constrained}: A builder makes hub investment decisions, considering the overall cost of constructing and operating the hub, such that the minimum-cost operation results in emissions at or below a given emissions target.
    \item \textbf{Framework 4: Tri-Level Regulator/Builder/Operator, Social Cost of Carbon}: A regulatory agency sets a Social Cost of Carbon (SCoC) rate, to reduce emissions to or below a certain target. The builder makes hub investment decisions considering the overall cost of constructing and operating the hub as well as the social cost of carbon emissions. The hub is independently operated based on a minimum cost-solution.
\end{itemize}

%\begin{figure}[h]
\begin{figure}
    \begin{center}
        \begin{tabular}{ r|c|c| }
        \multicolumn{1}{r}{}
         &  \multicolumn{2}{c}{builder \& operator} \\ \cline{2-3}
        \multicolumn{1}{r}{}
         &  \multicolumn{1}{c}{singular}
         & \multicolumn{1}{c}{distinct} \\
        \cline{2-3}
        carbon constrained & Framework 1 & Framework 3 \\
        \cline{2-3}
        carbon tax & Framework 2 & Framework 4 \\
        \cline{2-3}
        \end{tabular}
    \end{center}
    \caption{Taxonomy of low-carbon investment and design frameworks.}
    \label{frameworks}
\end{figure}

Each framework represents a different policy strategy that a regulator may have available for controlling carbon emissions. Each framework results in a distinct optimization formulation, with a distinct Pareto frontier of cost/emissions tradeoffs. Framework 1 represents the most efficient scenario, where a regulator can set an emissions limit and the builder and operator will work together to achieve it. Framework 2 represents the case where the regulator cannot dictate an emissions limit to the hub builder, but can set a carbon tax rate. Generally, higher tax rates are politically unpopular and increase the rate of tax evasion \cite{liu2013}. Therefore for the purposes of this framework, the regulator wants to minimize the tax rate (maximum likelihood of political feasibility), subject to the constraint that the chosen tax rate, if implemented, will result in meeting the emissions target. Framework 3 represents the case where the regulator sets an emissions limit and mandates that the builder chooses equipment such that the emissions constraint is robust to an operator that may decide to ignore emissions in favor of cost savings. Framework 4 is similar to Framework 3, except that the regulator cannot dictate an emissions limit to the hub operator, but can mandate that the builder consider a SCoC rate when deciding equipment investments. Each framework is further discussed in the following subsections.

\subsection{Common Framework Definitions}

Common between all framework formulations are the energy hub network constraints described in Section \ref{network_constraints}, as well as the definitions of the investment cost, operating cost, revenue, and emissions given in \eqref{invest_def}-\eqref{emission_def}, respectively. Eq. \eqref{invest_def} defines the total investment cost $C^{\text{invest}}$ in terms of the investment quantity variables multiplied by each variable's per-unit cost. The annual operating cost $C_{y}^{\text{operate}}$ depends on the grid power purchases and the fuel prices in each time period \eqref{operate_def}. The annual operating revenue $R_{y}^{\text{operate}}$ depends on the grid power exports and the feed-in prices in each time period \eqref{revenue_def}. Finally, the annual carbon emissions $E_{y}^{\text{operate}}$ depend on the grid power purchases and the carbon intensity of each fuel in each time period \eqref{emission_def}.

\begin{align}
    C^{\text{invest}} := & \sum_{m \in M} C_{m}^{\text{cap}} P_{m}^{\text{max}} + \sum_{g \in G^{\text{C}}} C_{g}^{\text{unit}} I_{g} \nonumber \\
    + & \sum_{g \in G^{\text{S}}} \big( C_{g}^{\text{power}} D_{g}^{\text{max}} + C_{g}^{\text{energy}} Q_{g}^{\text{max}} \big) \label{invest_def} \displaybreak[0] \\
    C_{y}^{\text{operate}} := 365 & \sum_{s \in S} \pi_{s} \sum_{t \in T} \sum_{m \in M} f_{m,s,t,y} P_{m,s,t,y} \Delta t \label{operate_def} \displaybreak[0] \\
    R_{y}^{\text{operate}} := 365 & \sum_{s \in S} \pi_{s} \sum_{t \in T} \sum_{m \in M} h_{m,s,t,y} r_{m,s,t,y} \Delta t \label{revenue_def} \displaybreak[0] \\
    E_{y}^{\text{operate}} := 365 & \sum_{s \in S} \pi_{s} \sum_{t \in T} \sum_{m \in M} e_{m,s,t,y} P_{m,s,t,y} \Delta t \label{emission_def}
\end{align}

\begin{figure*}[h]
    \centering
    \includegraphics[width=0.9\linewidth]{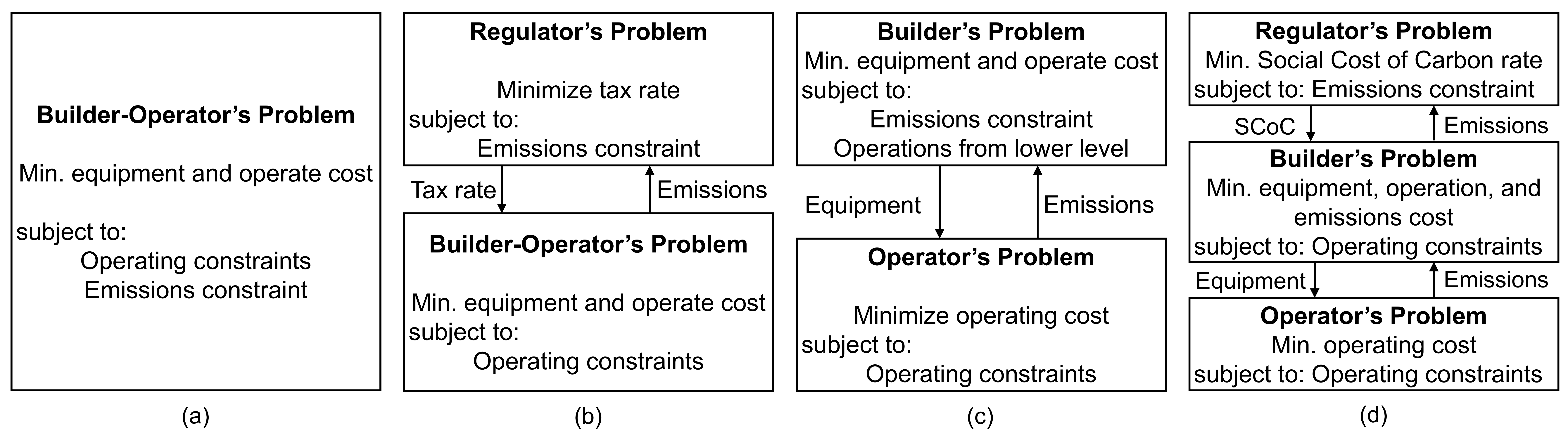}
    \caption{Low-carbon design frameworks: (a) Single Builder-Operator, (b) Bi-level Regulator/Builder-Operator, (c) Bi-level Builder/Operator, and (d) Tri-level Regulator/Builder/Operator.}
    \label{four_in_a_row}
\end{figure*}

\subsection{Framework 1: Single Builder-Operator}

The single builder-operator model features builder-operator coordination and an explicit carbon constraint, as shown in Fig. \ref{four_in_a_row}(a). The objective function is the net-present cost of building and operating the hub \eqref{sbo_obj}; Eqs. \eqref{sbo_obj}-\eqref{sbo_e_constraint} formalize this problem.

\begin{gather}
    \min C^{\text{invest}} + \sum_{y=1}^{Y} \frac{\big( C_{y}^{\text{operate}} - R_{y}^{\text{operate}} \big)}{\big( 1 + i \big)^{y}} \label{sbo_obj} \\
    \text{subject to:} \hspace{180pt} \nonumber \\
    \text{Network constraints: \eqref{converter_conservation}-\eqref{non_neg_flows}} \label{sbo_hub_constraints} \\
    E_{y}^{\text{operate}} \le E^{\text{max}} \hspace{8pt} \forall y \in Y \label{sbo_e_constraint}
\end{gather}

\subsection{Framework 2: Bi-Level, Regulator/Builder-Operator}

The bi-level regulator/builder-operator problem features builder-operator coordination and a carbon price set by an upper-level regulator, as shown in Fig. \ref{four_in_a_row}(b). Eqs. \eqref{tsbo_ul_obj}-\eqref{tsbo_hub_constraints} formalize this problem. In this framework, the regulator's objective is to find the minimum tax-rate $P^{\text{CO}_{2}}$ \eqref{tsbo_ul_obj} such that emissions are at or below a target \eqref{tsbo_e_constraint} when the builder-operator independently minimizes its net-present cost of building and operating the hub \eqref{tsbo_ll_obj}-\eqref{tsbo_hub_constraints} based on the tax-rate.

\begin{gather}
    \min P^{\text{CO}_{2}} \label{tsbo_ul_obj} \\
    \text{subject to:} \hspace{180pt} \nonumber \\
    E_{y}^{\text{operate}} \le E^{\text{max}} \hspace{8pt} \forall y \in Y \label{tsbo_e_constraint} \\
    E_{y}^{\text{operate}} \in \underset{\boldsymbol{I}, \boldsymbol{P}^{\text{max}}, \boldsymbol{D}^{\text{max}}, \boldsymbol{Q}^{\text{max}}, \boldsymbol{V}, \boldsymbol{Q}, \boldsymbol{r}}{\text{arg min}} \big \{ \hspace{118.9pt} \nonumber \\
    \hspace{40pt} C^{\text{invest}} + \sum_{y=1}^{Y} \frac{\big( C_{y}^{\text{operate}} - R_{y}^{\text{operate}} + T_{y}^{\text{operate}} \big)}{\big( 1 + i \big)^{y}} \label{tsbo_ll_obj} \\
    \text{subject to:} \hspace{100pt} \nonumber \\
    T_{y}^{\text{operate}} = \text{max} \left ( P^{\text{CO}_{2}} E_{y}^{\text{operate}}, 0 \right ) \hspace{8pt} \forall y \in Y \label{tsbo_tax_def} \\
    \text{Network constraints: \eqref{converter_conservation}-\eqref{non_neg_flows}} \label{tsbo_hub_constraints} \big \}
\end{gather}

The annual tax bill for carbon emissions $T_{y}^{\text{operate}}$ is constrained to be non-negative in \eqref{tsbo_tax_def}; otherwise, at high tax rates the lower-level problem can become unbounded. This can occur if there is at least one time period with negative marginal emissions of electricity, a condition that can be caused by transmission grid congestion.

\subsection{Framework 3: Bi-Level, Builder/Operator}

The bi-level builder/operator problem features no builder-operator coordination and an explicit carbon constraint, as shown in Fig. \ref{four_in_a_row}(c). Eqs. \eqref{bo_ul_obj}-\eqref{bo_hub_constraints} formalize this problem. The builder's objective is to minimize the overall construction and operating cost \eqref{bo_ul_obj} such that the emissions are at or below a target \eqref{bo_e_constraint} when the operator independently minimizes its operating cost \eqref{bo_ll_obj} subject to hub constraints \eqref{bo_hub_constraints}.
%The formulation of the builder/operator problem is given by:

\begin{gather}
    \min_{\boldsymbol{I}, \boldsymbol{P}^{\text{max}}, \boldsymbol{D}^{\text{max}}, \boldsymbol{Q}^{\text{max}}} C^{\text{invest}} + \sum_{y=1}^{Y} \frac{\big( C_{y}^{\text{operate}} - R_{y}^{\text{operate}} \big)}{\big( 1 + i \big)^{y}} \label{bo_ul_obj} \\
    \text{subject to:} \hspace{180pt} \nonumber \\
    \sum_{s \in S} \pi_s E_{s,y}^{\text{operate}} \le E^{\text{max}} \hspace{8pt} \forall y \in Y \label{bo_e_constraint} \\
    C_{s,y}^{\text{operate}}, E_{s,y}^{\text{operate}}, R_{s,y}^{\text{operate}} \in \underset{\boldsymbol{V}, \boldsymbol{Q}, \boldsymbol{r}}{\text{arg min}} \big \{ C_{s,y}^{\text{operate}} - R_{s,y}^{\text{operate}} \label{bo_ll_obj}
\end{gather}
\hspace{65pt} subject to:
\begin{gather}
    \text{Network constraints: \eqref{converter_conservation}-\eqref{non_neg_flows}} \big \} \hspace{8pt} \forall s \in S, y \in Y \label{bo_hub_constraints}
\end{gather}

\subsection{Framework 4: Tri-Level, Regulator/Builder/Operator}

The tri-level Regulator/Builder/Operator problem features no builder-operator coordination and a price to be set for the Social Cost of Carbon, $SCoC$ (\textyen/ton), as shown in Fig. \ref{four_in_a_row}(d). Eqs. \eqref{rbo_ul_obj}-\eqref{rbo_hub_constraints} formalize this problem. The regulator's objective is to minimize the $SCoC$ \eqref{rbo_ul_obj} such that the builder, when minimizing the combined cost of investment, operation, and carbon \eqref{rbo_ml_obj}, and the operator, minimizing its operating cost \eqref{rbo_ll_obj} subject to energy hub constraints \eqref{rbo_hub_constraints}, result in annual emissions below a given target \eqref{rbo_e_constraint}. Note that the $SCoC$ is not a cost that necessarily needs to be paid, as long as it is considered in the builder's objective function.

\vspace{-12pt}
\begin{gather}
    \min SCoC \label{rbo_ul_obj} \\
    \text{subject to:} \hspace{180pt} \nonumber \\
    \sum_{s \in S} \pi_s E_{s,y}^{\text{operate}} \le E^{\text{max}} \hspace{8pt} \forall y \in Y \label{rbo_e_constraint} \\
    E_{s,y}^{\text{operate}} \in \underset{\boldsymbol{I}, \boldsymbol{P}^{\text{max}}, \boldsymbol{D}^{\text{max}}, \boldsymbol{Q}^{\text{max}}}{\text{arg min}} \Big \{ \nonumber \\
    C^{\text{invest}} + \sum_{y=1}^{Y} \left [ \frac{\big( C_{y}^{\text{operate}} - R_{y}^{\text{operate}} \big)}{\big( 1 + i \big)^{y}} + SCoC \cdot E_{y}^{\text{operate}} \right ] \label{rbo_ml_obj} \\
    \text{subject to:} \hspace{80pt} \nonumber \\
    C_{s,y}^{\text{operate}}, E_{s,y}^{\text{operate}}, R_{s,y}^{\text{operate}} \in \underset{\boldsymbol{V}, \boldsymbol{Q}, \boldsymbol{r}}{\text{arg min}} \big \{ C_{s,y}^{\text{operate}}- R_{s,y}^{\text{operate}} \label{rbo_ll_obj} \\
    \text{subject to:} \nonumber \\
    \text{Network constraints: \eqref{converter_conservation}-\eqref{non_neg_flows}} \big \} \hspace{8pt} \forall s \in S, y \in Y \Big \} \label{rbo_hub_constraints}
\end{gather}

\section{Energy Hub Model Formulation} \label{network_constraints}

The energy hub operational model assumes a greenfield design, where the topology of the network is not predefined \cite{wang2017apen}, with enhancements to reduce dimensionality and include additional investment decisions. In the greenfield model in \cite{wang2017apen}, there is a branch between every combination of device output port and device input port which handle the same energy type. For that formulation, the number of branch flow variables grows with the number of devices at a rate of $\mathcal{O}(n^2)$.

\begin{figure}
    \centering
    \includegraphics[width=\linewidth]{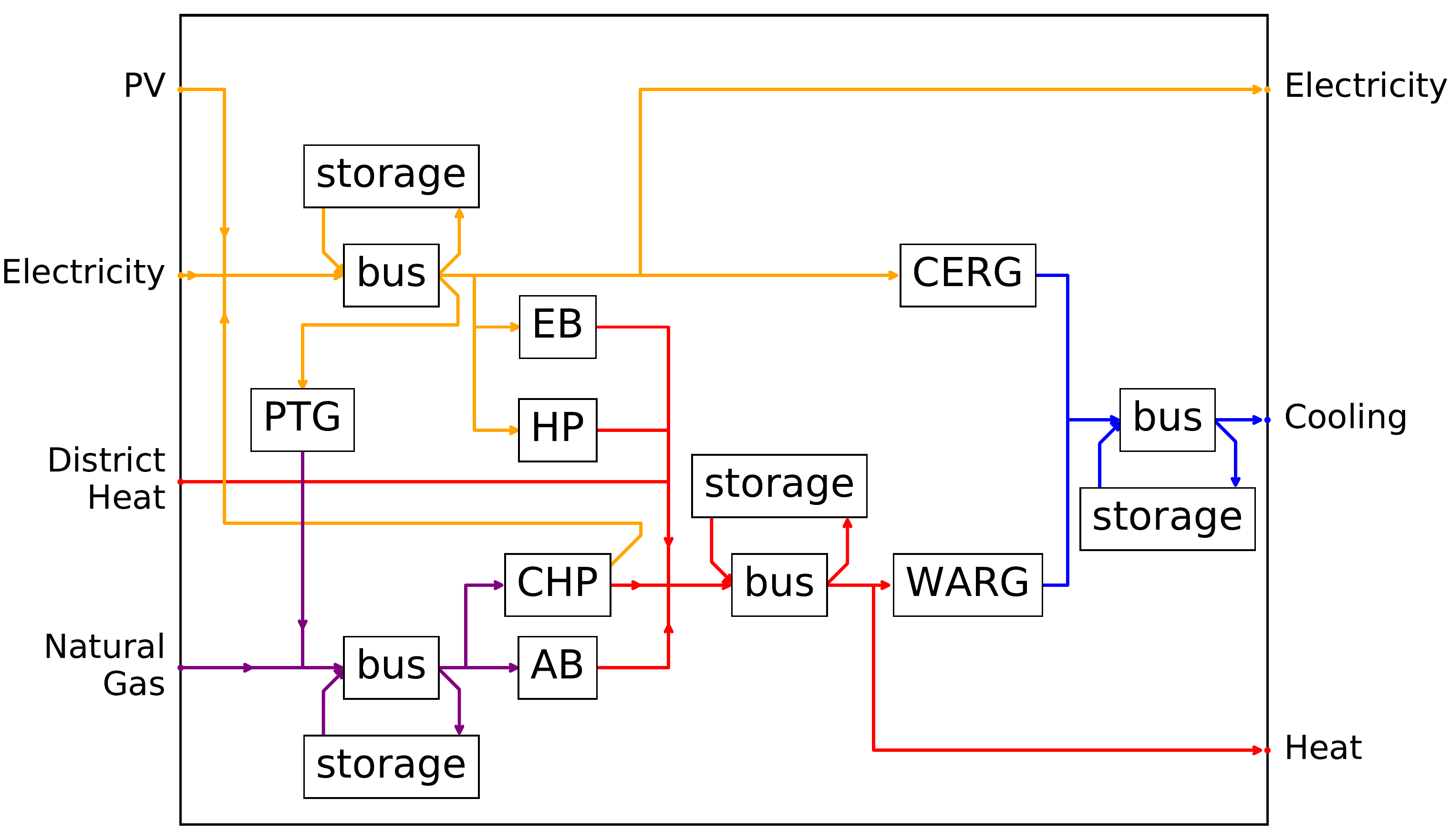}
    \caption{Network topology. Abbreviations: auxiliary boiler (AB), compression expansion refrigeration group (CERG), combined heat and power (CHP), electric boiler (EB), heat pump (HP), power-to-gas (PTG), water absorption refrigeration group (WARG).}
    \label{topology}
\end{figure}

By contrast, this formulation introduces the concept of \textit{energy buses}. Each energy bus is positioned such that any sources of this energy (i.e., from hub input ports and from converters and storage devices) flow directly into this bus, and any sinks of this energy (to hub output ports and to converters and storage devices) are fed directly from this bus. Therefore, the number of branch flow variables is reduced when compared to a formulation where every source can connect to every sink; $\mathcal{O}(n)$ vs. $\mathcal{O}(n^2)$. Fig. \ref{topology} shows the resulting network topology.

In addition to the converters in \cite{wang2017apen}, investment decisions include on-site distributed renewable generation, the provision of grid import/export capacity, storage for all energy types, and the inclusion of power-to-gas equipment \cite{vandewalle2015}. Since $\text{CO}_{2}$ is an input to the power-to-gas process, gas which is created using this process and burned is $\text{CO}_{2}\text{-neutral}$, though the electricity used during the conversion process may not be. For the considered devices, connecting all compatible input and output ports for electricity flows would require 23 branches, while introducing energy buses enables the same functionality with only 10 branches.

%\begin{figure}[h]

For a given energy hub, the operating decision variables are the branch flows $\boldsymbol{V}_{t}$, the state of charge for the storage devices $\boldsymbol{Q}_{t}$, and the grid export $\boldsymbol{r}_{t}$. Eqs. \eqref{converter_conservation}-\eqref{non_neg_flows} constrain the values of these decision variables. The dual variable associated with each constraint is shown in parentheses.

\vspace{-12pt}
\begin{gather}
    \boldsymbol{Z} \boldsymbol{V}_{t} = \boldsymbol{0} \hspace{8pt} \forall t \in T \hspace{8pt} (\alpha_{p,t}) \label{converter_conservation} \\
    Q_{g,t} = Q_{g,t-1} - \boldsymbol{J}_{g} \boldsymbol{A}_g \boldsymbol{V}_{t} \Delta t \hspace{8pt} \forall g \in G^{\text{S}}, t \in T \hspace{8pt} (\beta_{g,t}) \label{soc} \\
    0 \le Q_{g,t} \le Q_g^{\text{max}} \hspace{8pt} \forall g \in G^{\text{S}}, t \in T \hspace{8pt} (\underline{\gamma}_{g,t}, \overline{\gamma}_{g,t}) \label{soc_limit} \displaybreak[0] \\
    \boldsymbol{J}_{g} \boldsymbol{A}_g \boldsymbol{V}_{t} \le D_g^{\text{max}} I_{g} \hspace{8pt} \forall g \in G^{\text{C}}, t \in T \hspace{8pt} (\zeta_{g,t}) \label{converter_capacity} \displaybreak[0] \\
    - D_g^{\text{max}} \le \boldsymbol{J}_g \boldsymbol{A}_g \boldsymbol{V}_{t} \le D_g^{\text{max}} \hspace{8pt} \forall g \in G^{\text{S}}, t \in T \hspace{8pt} (\underline{\kappa}_{g,t}, \overline{\kappa}_{g,t}) \label{storage_power_capacity} \displaybreak[0] \\
    P_{m,t} := \boldsymbol{U} \boldsymbol{V}_{t} \le B_{m,t} P_m^{\text{max}} \hspace{8pt} \forall m \in M, t \in T \hspace{8pt} (\rho_{m,t}) \label{inputs_capacity} \displaybreak[0] \\
    \boldsymbol{W} \boldsymbol{V}_{t} = L_{m,t} + r_{m,t} \hspace{8pt} \forall m \in M, t \in T \hspace{8pt} (\mu_{m,t}) \label{outputs_def} \displaybreak[0] \\
    0 \le r_{m,t} \le P_{m}^{\text{max}} \hspace{8pt} \forall m \in M, t \in T \hspace{8pt} (\underline{\phi}_{m,t}, \overline{\phi}_{m,t}) \label{grid_sales_constraints} \displaybreak[0] \\
    0 \le K_{l} V_{l,t} \hspace{8pt} \forall t \in T, l \in L \hspace{8pt} (\sigma_{l,t}) \label{non_neg_flows} \displaybreak[0] \\
    \forall s \in S, y \in Y \nonumber
\end{gather}

Eq. \eqref{converter_conservation} expresses the conservation of power in all converters, storage devices, and energy buses using a connection and efficiency matrix $\boldsymbol{Z}$ \cite{wang2017apen}. The state of charge for these storage devices is tracked in \eqref{soc} and bounded by energy capacities in \eqref{soc_limit}. Power flows through each converter are constrained by \eqref{converter_capacity} and through each storage device by \eqref{storage_power_capacity}. Eq. \eqref{inputs_capacity} constrains hub input flows based on grid connection capacity and time-varying power availability, while \eqref{outputs_def} relates hub output flows to end-use power demand and power exports and \eqref{grid_sales_constraints} ensures that power exports do not exceed the grid connection capacity. Eq. \eqref{non_neg_flows} ensures that the directionality of branch flows is maintained, except for state of charge branches which are bi-directional.

Eqs. \eqref{converter_conservation}-\eqref{non_neg_flows} apply for each representative day, during each year. Therefore an additional two dimensions ($s$ and $y$) are added to variables $\boldsymbol{V}_{t}$, $\boldsymbol{Q}_{t}$, and $\boldsymbol{r}_{t}$ and to parameters $B_{m,t}$ and $L_{m,t}$; Eqs. \eqref{converter_conservation}-\eqref{non_neg_flows} then apply $\forall s \in S, y \in Y$.

\section{Multi-Level Reformulation Techniques} \label{bi_level_to_single_level}

Two techniques are used to solve the multi-level optimization problems presented in Section \ref{design_frameworks}. The first technique constructs a single-level equivalent MILP from the bi-level Builder/Operator problems of Frameworks 3 and 4. There are two steps to this technique: first a non-linear single-level equivalent is constructed (Section \ref{bi_level_to_single_level}.A), and then it is approximated by a MILP formulation (Section \ref{bi_level_to_single_level}.B).

The second technique solves the minimum carbon-price problem with an upper-level regulator and a lower-level non-convex builder-operator using the bisection method (Section \ref{bi_level_to_single_level}.C). This technique is used in Framework 2 and Framework 4, after a single-level equivalent is constructed.

For each framework, since the operation problems for each day are independent, the problem can be decomposed using Benders' Decomposition \cite{benders1962}: the investment variables are solved in the master problem and each day's operational variables are solved in a subproblem.

\subsection{Constructing a Single-Level Equivalent for the Distinct Builder and Operator Case}

An independent cost-minimization problem for the operator must be solved for each representative day $s$ of each representative year $y$. For investment problem formulations where the builder and the operator are distinct, the builder must anticipate the emissions resulting from these cost-minimization decisions to make its investment decisions; this ensures that the emissions constraints can be satisfied or the social cost of GHG emissions can be appropriately considered in the objective. Since the operator's cost minimization problem is linear and therefore convex, the strong duality theorem can be used to constrain the variables in the operator's problem (the most relevant of which are the resulting cost and emissions) to only the set of cost-minimizing values \cite{boyd2004}.

The dual problem is defined in terms of the dual variables of the constraints \eqref{converter_conservation}-\eqref{non_neg_flows} and the parameters of the original primal problem. The dual objective is given in \eqref{dll_obj} with feasibility constraints corresponding to each primal variable given in \eqref{dll_q_eq} for $\boldsymbol{Q}_t$, \eqref{dll_r_eq} for $\boldsymbol{r}_t$, and \eqref{dll_v_eq} for $\boldsymbol{V}_t$, with dual variable domains given in \eqref{dll_domain}.

\begin{align}
    & \max_{\alpha, \beta, \gamma, \zeta, \kappa, \rho, \mu, \phi, \sigma} \hat{C}^{\text{DLL}} := 
    \sum_{t \in T} \bigg \{ \sum_{g \in G^{\text{C}}} \Big ( - D_{g} I_{g} \zeta_{g,t} \bigg ) \nonumber \\
    & + \sum_{g \in G^{\text{S}}} \left (
    - Q_{g}^{\text{max}} \overline{\gamma}_{g,t}
    - D_g^{\text{max}} \left ( \underline{\kappa}_{g,t} + \overline{\kappa}_{g,t} \right ) \right ) \nonumber \\
    & + \sum_{m \in M} \left ( - L_{m,t} \mu_{m,t} - P_{m}^{\text{max}} \left ( B_{m,t} \rho_{m,t} + \overline{\phi}_{m,t} \right ) \right ) \bigg \} \label{dll_obj}
\end{align}

subject to:
\vspace{-6pt}
\begin{align}
    \beta_{g,t} &- \beta_{g,t+1} + \overline{\gamma}_{g,t} - \underline{\gamma}_{g,t} = 0
    \hspace{12pt} \forall g \in G^{\text{S}},  t \in T
    \label{dll_q_eq} \displaybreak[0] \\
    \overline{\phi}_{m,t} &- \underline{\phi}_{m,t} - h_{m,t} \Delta t - \mu_{m,t} = 0
    \hspace{6pt} \forall m \in M, t \in T
    \label{dll_r_eq} \displaybreak[0] \\
    \sum_{m \in M} & \bigg [ U_{m,l} \Big ( f_{m,t} \Delta t + \rho_{m,t} \Big ) + W_{m,l} \mu_{m,t} \bigg ] \nonumber \\
    + \sum_{p \in P} & \bigg \{
    \sum_{g \in G^{\text{S}}} \Big [ J_{g,p} A_{p,l} \big( \Delta t \beta_{g,t} + \overline{\kappa}_{g,t} - \underline{\kappa}_{g,t} \big ) \Big ] \nonumber \\
    + \sum_{g \in G^{\text{C}}} & J_{g,p} A_{p,l} \zeta_{g,t} \bigg \}
    + \sum_{p \in P^{\text{out}}} \Big (Z_{l,p} \alpha_{p,t} \Big ) - K_{l} \sigma_{l,t} = 0 \nonumber \\
    & \hspace{135pt} \forall l \in L, t \in T \label{dll_v_eq} \displaybreak[0] \\
    & \hspace{20pt} \underline{\gamma}, \overline{\gamma}, \zeta, \underline{\kappa}, \overline{\kappa}, \rho, \underline{\phi}, \overline{\phi}, \sigma \ge 0 \label{dll_domain}
\end{align}

Therefore, the power flows within the energy hub on a given day can be constrained to the values that would result from a hub operator's cost minimization using the primal feasibility constraints \eqref{primal_constraints}, the dual feasibility constraints \eqref{dual_constraints}, and the strong duality constraint \eqref{strong_duality}, and the problem can simultaneously consider the upper-level objective and constraints of the regulator.

\begin{gather}
    \text{Eqs. \eqref{converter_conservation}-\eqref{non_neg_flows} \hspace{12pt} (PLL feasibility)} \label{primal_constraints} \\
    \text{Eqs. \eqref{dll_q_eq}-\eqref{dll_domain} \hspace{12pt} (DLL feasibility)} \label{dual_constraints} \\
    \hat{C}^{\text{PLL}} := C^{\text{operate}} - R^{\text{operate}} = \hat{C}^{\text{DLL}} \hspace{12pt} \text{(Strong duality)} \label{strong_duality} %\\
    %\forall s \in S, y \in Y \nonumber
\end{gather}

\subsection{MILP Approximation of Non-Linear Strong Duality Constraint} \label{milp_approx}

Although the hub operator's cost-minimization problem, its dual problem, and the strong duality constraint are all linear in terms of the lower-level primal and dual variables, the strong duality constraints are non-linear when the equipment capacities are included as decision variables. To avoid the difficulty of solving a Mixed Integer Non-Linear Problem (MINLP), these equipment capacity variables can be discretized in order to convert the problem to a more tractable MILP version.

The continuous variables $P_{g}^{\text{max}}$, $Q_{g}^{\text{max}}$, $D_{g}^{\text{max}}$, and the integer variables $I_{g}^{\text{max}}$ are approximated by a series of binary variables using \eqref{input_cap_sum}-\eqref{converter_invest_sum}. Rather than having each binary variable represent a single unit of capacity, a binary counting approach is used, where each binary variable represents one digit of a binary representation of an integer number, as shown in Fig. \ref{binary_representation}. In this way, a relative step size of $1/2^{n}$ is possible with only $n$ binary variables for each continuous variable (e.g., a step size of 0.1\% for 10 binaries). This method implicitly creates bounds on the approximated continuous variables, since each approximated value can only range from the value represented by $\{0,0,\dots,0\}$ to the value represented by $\{1,1,\dots,1\}$. Therefore, the values of the step sizes must be chosen carefully.

\begin{gather}
    %Summed binary variables
    P_{g}^{\text{max}} \approx \sum_{n=0}^{N} 2^{n} x_{m,n}^{\text{a}} \Delta P_{g} \hspace{8pt} \forall m \in M \displaybreak[0]
    \label{input_cap_sum} \\
    Q_{g}^{\text{max}} \approx \sum_{n=0}^{N} 2^{n} x_{g,n}^{\text{b}} \Delta Q_{g} \hspace{8pt} \forall g \in G^{\text{S}} \displaybreak[0]
    \label{storage_energy_sum} \\
    D_{g}^{\text{max}} \approx \sum_{n=0}^{N} 2^{n} x_{g,n}^{\text{c}} \Delta D_{g} \hspace{8pt} \forall g \in G^{\text{S}} \displaybreak[0]
    \label{storage_power_sum} \\
    I_{g}^{\text{max}} \approx \sum_{n=0}^{N} 2^{n} x_{g,n}^{\text{d}} \hspace{8pt} \forall g \in G^{\text{C}} \displaybreak[0]
    \label{converter_invest_sum}
    %\displaybreak[0] \\
\end{gather}

\begin{figure}[h]
    \begin{center}
        \begin{tabular}{|r|c|c|c|c|c|c|c|c|c|c|}
            \hline
            $681=$
            & \multicolumn{1}{c}{$2^9$}
            & \multicolumn{1}{c}{$+$}
            & \multicolumn{1}{c}{$2^7$}
            & \multicolumn{1}{c}{$+$}
            & \multicolumn{1}{c}{$2^5$}
            & \multicolumn{1}{c}{$+$}
            & \multicolumn{1}{c}{$2^3$}
            & \multicolumn{1}{c}{}
            & \multicolumn{1}{c}{$+$}
            & $2^0$
            \\ \hline
            $=$ & 1 & 0 & 1 & 0 & 1 & 0 & 1 & 0 & 0 & 1
            \\ \hline
        \end{tabular}
    \end{center}
    \caption{Binary representation of 681. 10 bits can represent 1024 values.}
    \label{binary_representation}
\end{figure}

Using \eqref{input_cap_sum}-\eqref{converter_invest_sum}, the non-linearity in the strong duality constraint is therefore reduced to only products of binary and continuous variables. These non-linear products are approximated using the big-M technique \cite{floudas1995}, resulting in strong duality constraints that are MILP rather than MINLP.
For computational performance, it is standard practice to choose the values of M such that they are as small as possible while still allowing for the full range of values for the continuous variables \cite{cplex}. Since the dual variables for the network constraints are related to the price of input energy flows, scaling energy prices (and therefore big-M values) may improve performance.

\subsection{Determining a Minimum Carbon Tax using the Bisection Method}

Incorporating a carbon price in the lower-level builder-operator problem is an instance of using a linear weighting factor to optimize trade-offs between multiple competing objectives. In this case, these objectives are the operating emissions and the sum of investment and operating costs. By varying the carbon price ($P^{\text{CO}_{2}}$ or $SCoC$), the Pareto frontier is sampled at points which lie on the convex hull of the cost/emissions solution space \cite{deb2016}. Since the lower-level builder-operator problem is non-convex due to the inclusion of binary variables, the minimum-carbon-price solution which satisfies an emissions constraint may result in a total investment and operating cost that is higher than could be found using an $\epsilon$-constraint method (Frameworks 1 \& 3). However, the resulting emissions would also be lower.

Since the upper-level price-setting problem minimizes a single continuous variable (the tax rate), and the resulting emissions from the lower-level problem are monotonically non-increasing with respect to an increasing tax rate, the problem can be solved to within a specified tolerance by the bisection method, as shown in \cite{olsen2018}.

\subsection{Computational Complexity}

Since the performance of MILP solvers depends on the number of variables and constraints in a model, the number of variables of each type are listen in Table \ref{computational_complexity}, where set names are used to represent the number of elements in the set. These expressions can be used to balance the accuracy of different model elements when computation time and/or available memory are an issue.

\begin{table}[h]
    \caption{Computational Complexity} 
    \begin{tabular}{ 
        			>{\centering\arraybackslash}p{0.35\linewidth}
        			>{\centering\arraybackslash}p{0.55\linewidth}
            		}
        \hline \hline
        \# of integer variables & $M + 2G^{\text{S}} + G^{\text{C}}$
        \\ \hline
        \# of binary variables &
        $N^{\text{A}} M
        + G^{\text{S}} \left ( N^{\text{B}} + N^{\text{C}} \right )
        + N^{\text{D}} G^{\text{C}} $
        \\ \hline
        \# of continuous variables &
        $STY \big [
        2 L + 6 G^{\text{S}} + 5 M + P^{\text{out}} + G^{\text{C}}$   %PLL/DLL
        $+ 2 M N^{\text{A}}
        + 2 G^{\text{S}} \left ( N^{\text{B}} + N^{\text{C}} \right )
        + G^{\text{C}} N^{\text{D}}
        \big ]$
        \\ \hline
        \# of constraints &
        $STY \big [
        2 L + 6 G^{\text{S}} + 5 M + P^{\text{out}} + G^{\text{C}}$
        $+ 6 M N^{\text{A}}
        + 6 G^{\text{S}} \left ( N^{\text{B}} + N^{\text{C}} \right )
        + 3 G^{\text{C}} N^{\text{D}}
        \big ]$
    \end{tabular}
    \label{computational_complexity}
\end{table}

\section{Case Study}

\subsection{Parameters}

This case study considers the construction of a new subsidiary administrative center in the Tongzhou region of Beijing. Fuel prices and converter parameters are taken from \cite{wang2017apen}, with electricity being bought-back at 85\% of the off-peak price. A power-to-gas converter is added, with an assumed capacity of 100 MW at a cost of 40 M\textyen\hspace{0.1em} \cite{desaintjean2014}. Increased prices of electricity during peak periods are in force from June 1st to September 31st. District heating is assumed to be available  only from November 15th to March 15th each year. Heating, cooling, and electricity demand patterns as a function of outdoor temperature are adapted from \cite{zhang2016}. Solar generation profiles and temperatures for Beijing are obtained from NREL's Typical Meteorological Year dataset \cite{tmy_data}. The optimization horizon is 20 years, with fuel prices increasing at 2\% per year, energy demands at 4\%, and a 10\% discount rate. All investments are made at year zero.

Input capacity is assumed to cost 100 k\textyen/MW for each fuel, and PV capacity is assumed to cost 5 M\textyen/MW. Electricity storage capacity is assumed to cost \$1000/kW and \$50/kWh \cite{pandzic2015}. Thermal storage capacity is assumed to cost \$25/kW and \$25/kWh \cite{hasnain1998}. Bulk LNG storage capacity costs approximately \$250/MWh \cite{aspelund2009} and \$20,000/MW \cite{songhurst2014} (converted from tons, and tons/year). Since a campus-scale LNG facility is smaller than the export-scale facilities in \cite{aspelund2009, songhurst2014}, specific costs are assumed to be greater by a factor of 10. All calculations are conducted in \textyen, at a rate of 6.6 \textyen/\$.

The carbon intensity of the district heating system is assumed to be 0.3 t/MWh \cite{werner2017}. The carbon intensity of natural gas is taken as 0.181 t/MWh, based on the chemical composition of methane. The carbon intensity of grid power is time-varying, based on the characteristics of the marginal generators \cite{graffzivin2014}. For this case study, a range of marginal emission rates is obtained by running unit commitment problems for a modified ISO-NE test system \cite{krishnamurthy2016}, with wind generation providing 15\% of the annual energy consumption.

Since solving the planning problem while modeling operations for each day of the year entails an excessively high amount of computation, a subset of representative days are selected using a modified $k$-means algorithm, where discrete variables are preserved and distances from clusters are evaluated based on the Z-score for each time-varying parameter. For this case study, a $k$ value of 10 was chosen because it appears to balance computational complexity with descriptiveness, with at least one time period of negative marginal emissions. Short-term marginal emissions can be zero due to renewable spillage, or even negative due to transmission congestion.

\subsection{Solution Methods}

The model is implemented in GAMS 25.0 \cite{gams} and solved using CPLEX 12.8 \cite{cplex} on machines with at least 16 cores, each running at at least 2 GHz, with at least 64 GB of RAM.

Using this combination of hardware and software, progress on closing the optimality gap can be slow for the distinct builder/operator frameworks (Frameworks 3 and 4), as these problems are in general NP-hard \cite{dempe2015}. Progress of branch and cut (B\&C) solvers can generally be improved by providing the solver with a feasible warm start, to provide an upper limit for pruning branches and as a starting point for relaxation induced neighbor search heuristics \cite{danna2005}. A feasible solution for a problem with a given $E^{\text{max}}$ can be used as a warm start for a problem with any greater value of $E^{\text{max}}$ (a relatively relaxed problem).

Due to the binary representation of integer variables described in Section \ref{milp_approx}, solutions representing `adjacent' integer values can be very different in terms of their binary components, and vice versa, which may hinder the effectiveness of neighbor-search heuristics built into B\&C solvers. Therefore, a branch-and-cut-and-heuristic \cite{gams_bch} process is implemented, in which incumbent solutions are periodically output to an accompanying heuristic which attempts to find a better feasible candidate to return to the B\&C solver.

For a candidate solution with $E^{\text{operate}} < E^{\text{max}}$, the set of adjacent solutions are first enumerated by perturbing the value of integer investment variables by one capacity step and then calculating the binary variable representation \textit{a priori}. These candidate adjacent solutions (`neighbors') are then evaluated for network feasibility, emissions-constraint feasibility, and cost reduction by running the networks constraints \eqref{converter_conservation}-\eqref{non_neg_flows} with investment variables fixed (an LP). If there are no adjacent solutions which meet the network feasibility constraint and emissions constraint at lower-cost, the heuristic halts. Otherwise, the feasible and lower-cost solutions are ranked in terms of their cost reduction per emissions-increase (neighbors with greater cost-reduction per emissions-increase are deemed `better'). From this subset of neighbors, the `best' candidate is chosen and the heuristic is repeated using that candidate as a starting point until a candidate is found for which there are no neighbors which meet the emissions and network feasibility constraints at lower cost. The lowest-cost emissions-feasible candidate is then returned to the B\&C solver as a new incumbent solution.

The computational complexity of these neighbor searches depends on the definition of `adjacent' solutions. If only one investment variable at a time is perturbed, the required number of LP solves will grow with $\mathcal{O}(n \cdot \Delta E)$, where $n$ is the number of integer investment variables ($|M|+|G|$) and $\Delta E$ is the difference between the emissions at the current candidate and the emissions target. If instead, neighbors are enumerated by perturbing two investment variables simultaneously (\textit{e.g.} the quantity of a particular type of converter is reduced while the capacity of a storage device is increased), the number of LP solves will grow with $\mathcal{O}(n^2 \cdot \Delta E)$; more neighbors are evaluated per heuristic iteration, providing potentially better performance per iteration at the expense of increased time per iteration.

\subsection{Sampling the Pareto Frontier}

In order to sample the Pareto frontier for the frameworks with direct emissions constraints (\textit{i.e.} Frameworks 1 and 3), the minimum cost emissions-unconstrained solution can first be obtained to determine the baseline emissions, and then a set of emissions-limits can be calculated for a specified resolution (e.g. 100 points for 1\% $E^{\text{max}}$ resolution); these independent emissions-limit problems can then be solved in parallel.

Information from `nearby' solutions can be used to improve knowledge about the optimality gap of the best known solution. For a given emissions target problem (A), if a solution for a tighter emissions-limit problem (B) is found with a lower cost than the currently best-known solution for problem (A), the cheaper solution from (B) be substituted for the currently best-known solution for (A), since a solution for a tighter problem (B) is always valid for a relaxed problem (A). Conversely, for a given emissions target problem (C), if a better lower bound is found for a more relaxed emissions-limit problem (D) is found, the lower bound from (D) can be substituted for the lower bound in (C), since the true optimum for a tightened problem (C) can only be greater than or equal to that of the relaxed problem (D).

In order to sample the Pareto frontier of the weighted-sum frameworks (\textit{i.e.} Frameworks 2 and 4), the optimization problem can be solved for a range of carbon prices, \textit{i.e.} the carbon tax rate in Framework 2 or the \textit{SCoC} rate for Framework 4. For a given carbon price, the best warm start from the set of previously found solutions can be evaluated by calculating the total cost (investment + operation + total carbon penalty) for the set of currently-known solutions, and providing the solver with the lowest total cost solution as a warm start.

Several other solution techniques for bi-objective mixed-integer optimization problems also exist; \cite{soylu2016} presents one based on the $\epsilon$-\textit{Tabu}-constraint method and also provides a review of several others.

\subsection{Computational Performance}

\begin{figure}
    \centering
    \includegraphics[width=\linewidth]{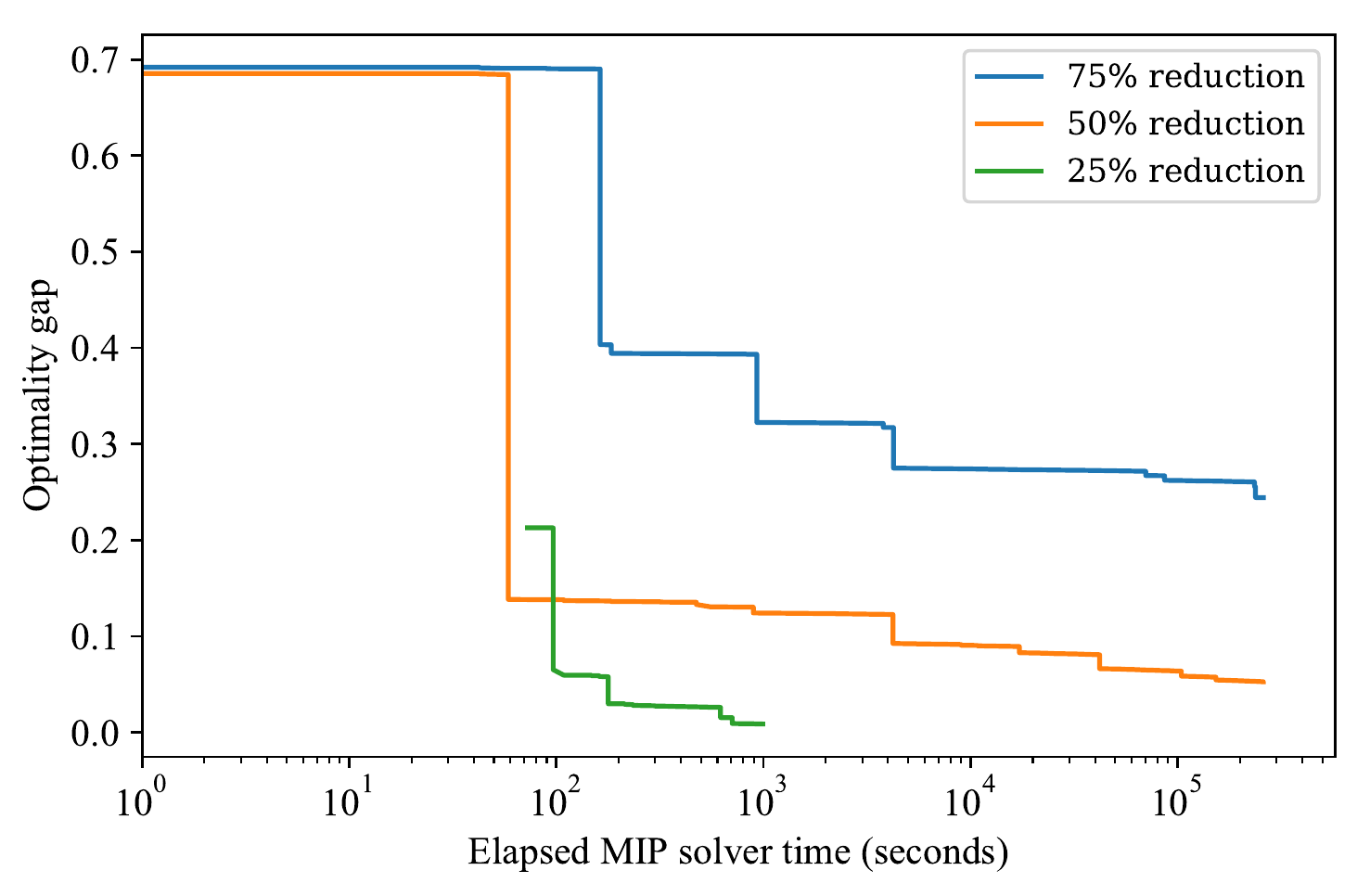}
    \caption{Optimality gap trajectories for selected emissions targets, Framework 3.}
    \label{trajectories}
\end{figure}

Frameworks 1 and 2 are solved to optimality for all emissions targets. When solving Framework 3, the performance of the MILP solver with the neighbor-search heuristic varies depends on the specified emissions target. For low emissions reductions targets, the problem can be solved to optimality; however, for moderate to aggressive targets CPLEX is found to stall in its progress in closing the optimality gap. For emissions reductions targets of up to 25\%, CPLEX can solve to within 1\% optimality gap. For emissions reductions targets of up to 50\%, CPLEX can solve to within 6\% optimality gap. The worst performance was found at an emissions reduction target of 72\%: progress stalls with an optimality gap of 17\%. Sample optimality gap trajectories are shown in Figure \ref{trajectories}.  Since the solution method for Framework 4 involves iterative solves of Framework 3, it suffers from the same computational challenges. To reduce the optimality gap further, the complexity of the case study must be reduced (see Table \ref{computational_complexity}), more processing power must be applied to the model, or more time must be allowed for solver convergence (the latter two of which are typically available to policy makers).

\section{Results}

\begin{figure}
    \centering
    \includegraphics[width=\linewidth]{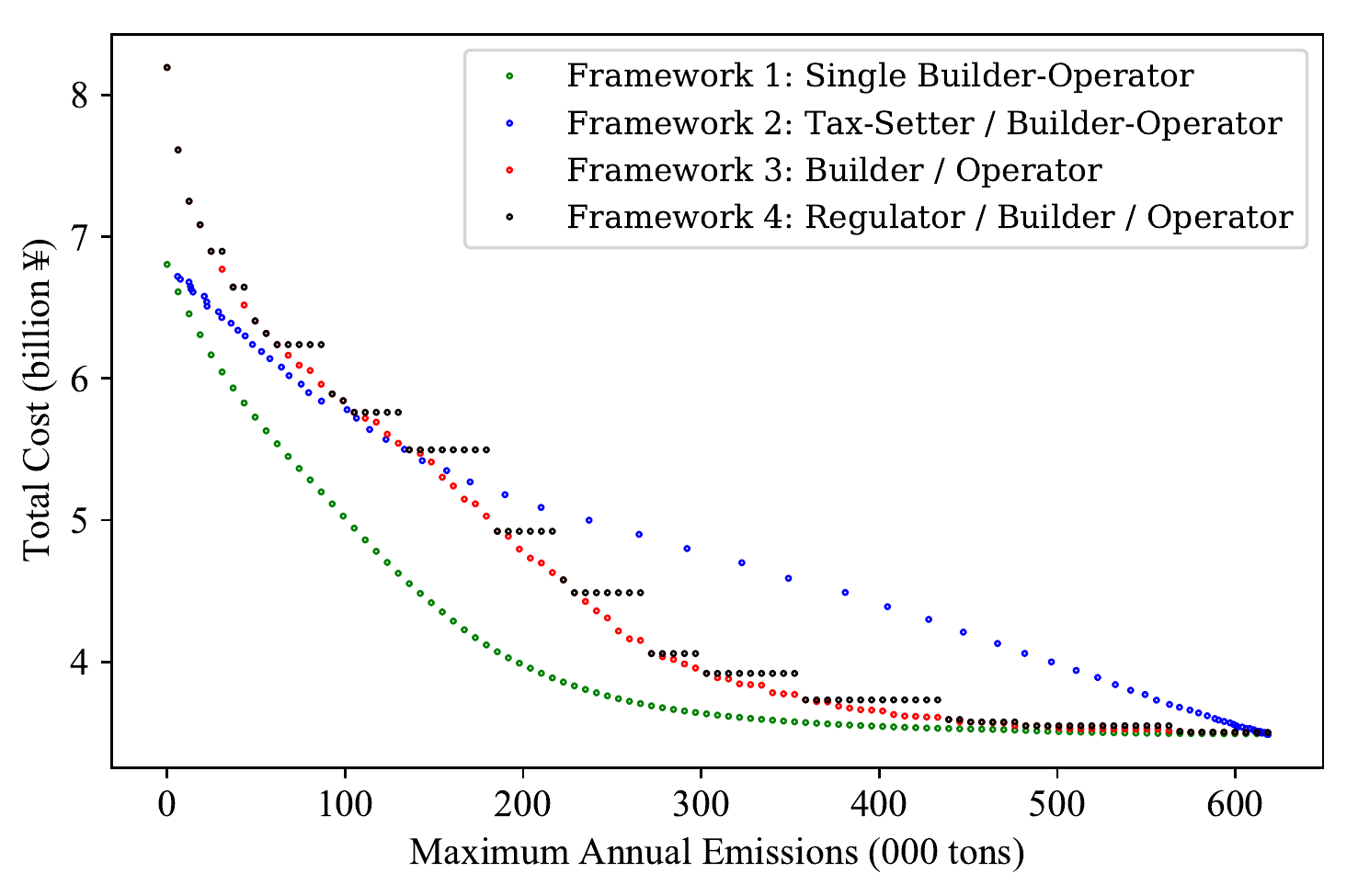}
    \caption{Comparison of cost/emissions Pareto frontiers for Frameworks 1-4.}
    \label{compare_methods}
\end{figure}

Fig. \ref{compare_methods} shows the total costs of the optimal energy hubs chosen by Frameworks 1-4 as a function of emissions target. Framework 1 results in the cheapest operating costs for a given emissions target, because the investment and operational variables are optimized simultaneously and there is no cost for emissions. Framework 2 can be significantly more expensive than Framework 1 for mid-range emissions reduction targets, but at very high emissions reductions targets this gap approaches zero as the high tax rate incentivizes the lower-level builder-operator to make investment and operating decisions with very low (or zero) carbon emissions. 
Framework 3 is not significantly more expensive than Framework 1 at modest emissions-reduction targets, but at aggressive emissions-reduction targets it is more expensive than both Framework 1 or Framework 2. In order to ensure that emissions decided by the lower-level operator do not exceed the target, the energy hub infrastructure must be `overbuilt': built to be able to satisfy end-use demands while giving the operator few to no opportunities to use cheap but carbon-intensive energy sources.
For a given emissions target, Framework 4 is as expensive or more expensive than Framework 3, since all of the constraints of the distinct Builder/Operator formulation still exist, but the upper-level regulator can only influence the investment decisions indirectly using the $SCoC$. This also creates large gaps between adjacent solutions. For example, no solutions are found with maximum annual emissions between 136 and 186 kilo-tons per year; any emissions target in this range must be met by imposing a high enough $SCoC$ to result in the (more expensive) 136 kilo-ton solution.

\subsection{Results for Framework 1: Single Builder-Operator}

In the single builder-operator framework, the total cost of building and operating the hub grows from 3.49 B\textyen\hspace{0.1em} when emissions are unconstrained (the `base case') to 6.80 B\textyen\hspace{0.1em} for a 100\% emissions reduction target. Cost components for several milestone emissions targets are shown in Table \ref{sbo_costs}. Much of the change in total cost occurs with aggressive emissions reductions targets, and consequently moderate emissions reductions targets are possible with only modest increases in total cost.

\begin{table}[h]
    \caption{Framework 1: Costs of solutions as a function of emissions targets} 
    \begin{tabular}{ 
        			>{\centering\arraybackslash}p{0.25\linewidth}
        			>{\centering\arraybackslash}p{0.15\linewidth}
        			>{\centering\arraybackslash}p{0.2\linewidth}
        			>{\centering\arraybackslash}p{0.2\linewidth}
            		}
        \hline \hline
        Emissions reduction target & Total cost (B\textyen) & Investment cost (B\textyen) & Net operational cost (B\textyen) \\ \hline
        None & 3.49 & 1.64 & 1.85 \\
        25\% & 3.52 & 1.62 & 1.90 \\
        50\% & 3.62 & 2.13 & 1.50 \\
        75\% & 4.35 & 4.12 & 0.23 \\
        100\% & 6.80 & 8.39 & -1.59
    \end{tabular}
    \label{sbo_costs}
\end{table}

\begin{figure}
    \centering
    \includegraphics[width=\linewidth]{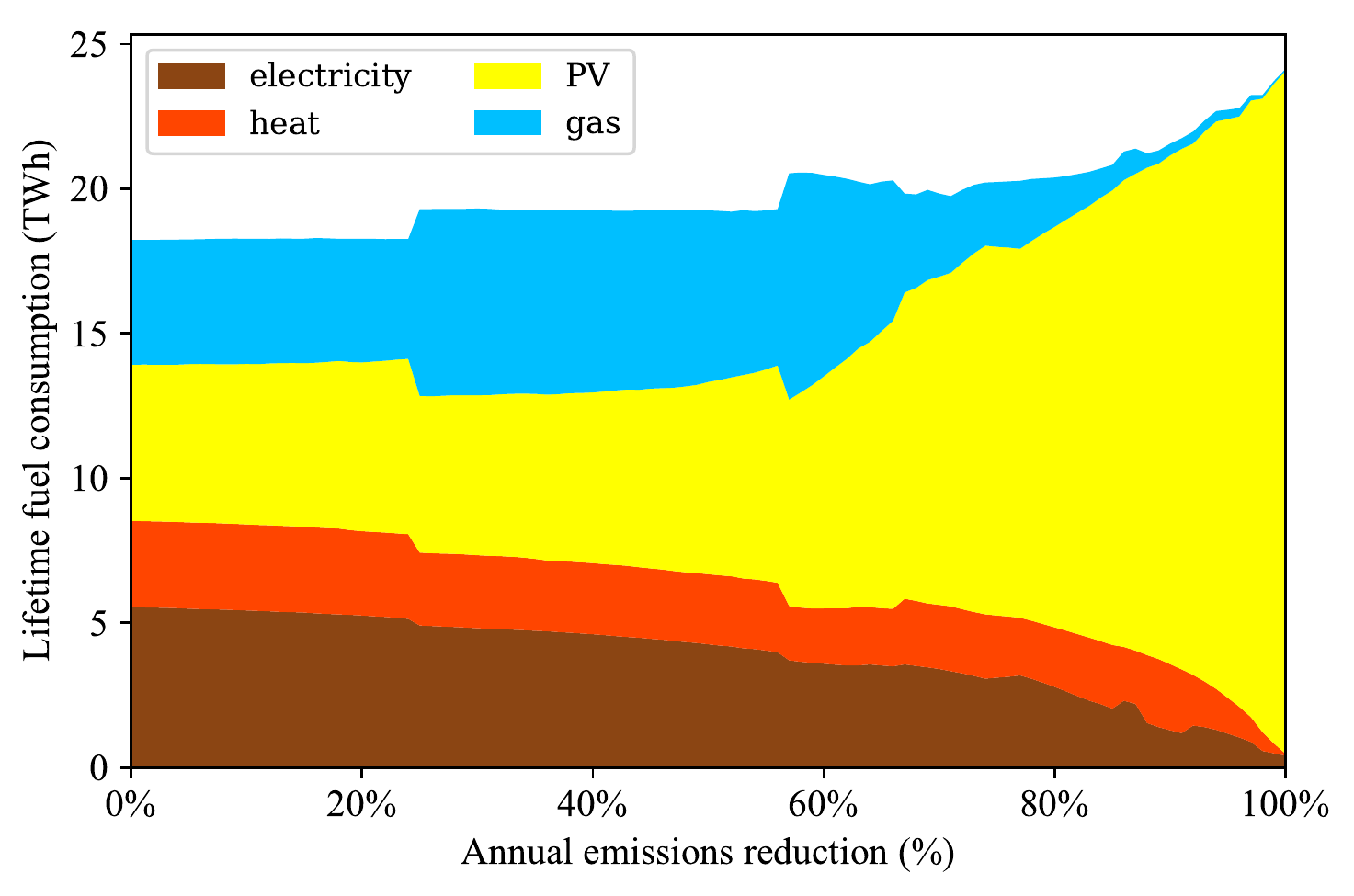}
    \caption{Fuel consumption as a function of emissions targets.}
    \label{fuel_consumption}
\end{figure}

Fig. \ref{fuel_consumption} shows the lifetime quantity of fuels flowing into the energy hub as a function of the emissions target. As the emission reduction target increases, the share of energy provided by electricity and district heat steadily declines. The changes in the costs components and in the fuel mix are significantly smaller in the 0\% to 50\% target range than in the 50\% to 100\% range. For moderate emissions reductions targets (up to approximately 65\%) gas consumption increases. On the other hand, past this point gas consumption drops quickly, to virtually zero in the zero-net-emissions case. PV generation increases by only 23\% between the base case and a 50\% target reduction, but grows by 337\% for a 100\% target reduction.

\begin{figure}
    \centering
    \includegraphics[width=\linewidth]{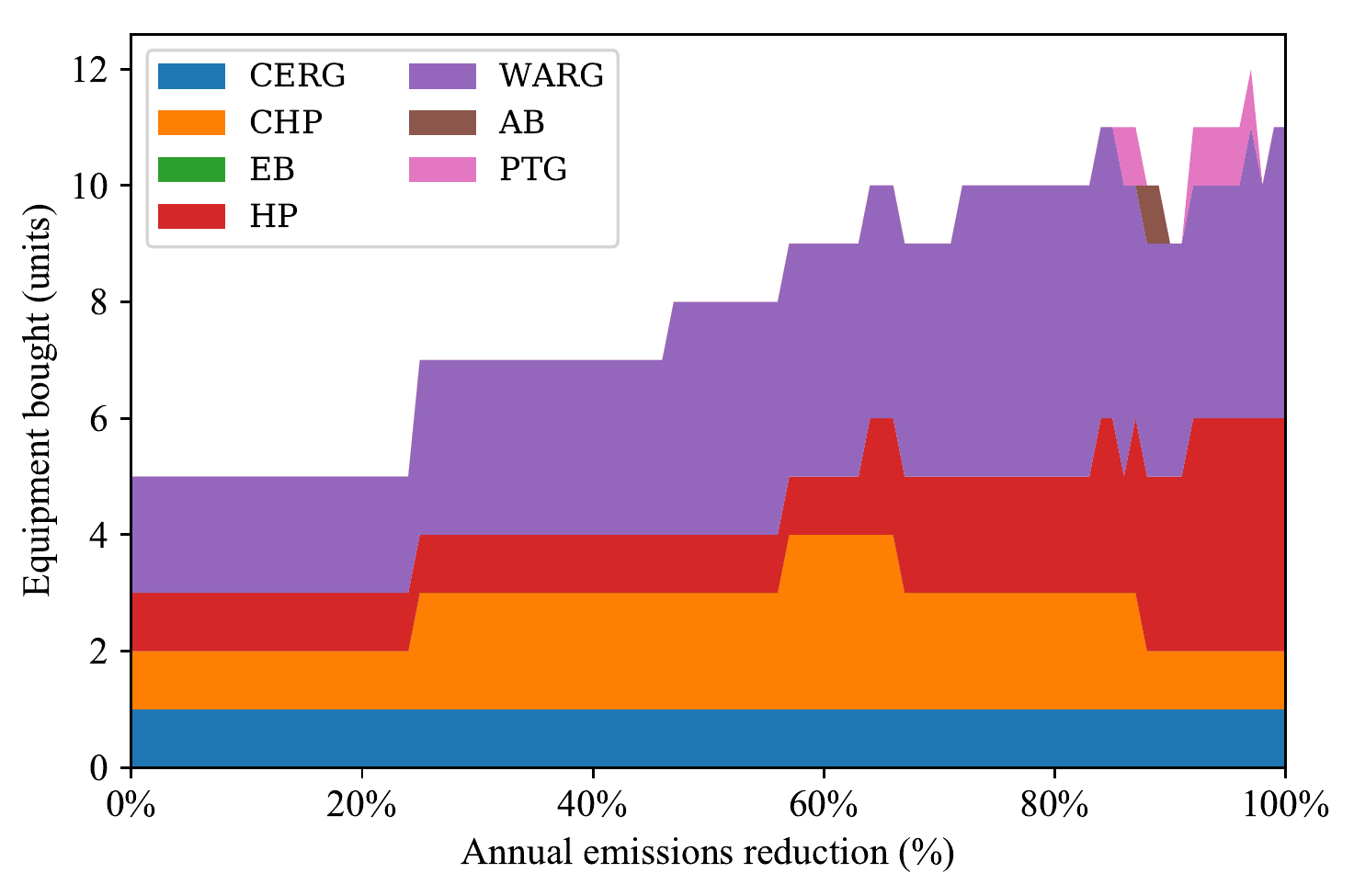}
    \caption{Equipment purchased as a function of desired emissions reduction. See Fig. \ref{topology} for abbreviations.}
    \label{equipment_purchased}
\end{figure}

Fig. \ref{equipment_purchased} shows how the choice of equipment varies as a function of the emissions reduction target. No matter the desired emissions reduction target, the optimal investment decisions include at least one 40 MW compression-expansion refrigeration group (CERG), one 100 MW combined-heat-and-power unit (CHP), one 40 MW heat pump (HP), and two 20 MW water-absorption refrigeration groups (WARG). The investment in CHP peaks at 3 units at emission reduction targets of 55-65\%, as the emissions created by burning gas are less than the emissions due to importing electricity and heat from the grid. Past this peak, the use of CHP declines and is replaced by electricity from PV generation and heat from heat pumps. The combination of HPs and WARGs is less efficient at converting electricity into cooling than CERGs, but the conversion to heat allows the use of intermediate thermal storage which is significantly cheaper than electricity storage. Storage allows end-use demands to be met with lower investments in input capacity and conversion equipment. Electric boilers are never chosen, and auxiliary boilers and power-to-gas units are only chosen for a few very aggressive emissions reduction targets.

\subsection{Results for Framework 2: Regulator / Builder-Operator}

\begin{figure}
    \centering
    \includegraphics[width=\linewidth]{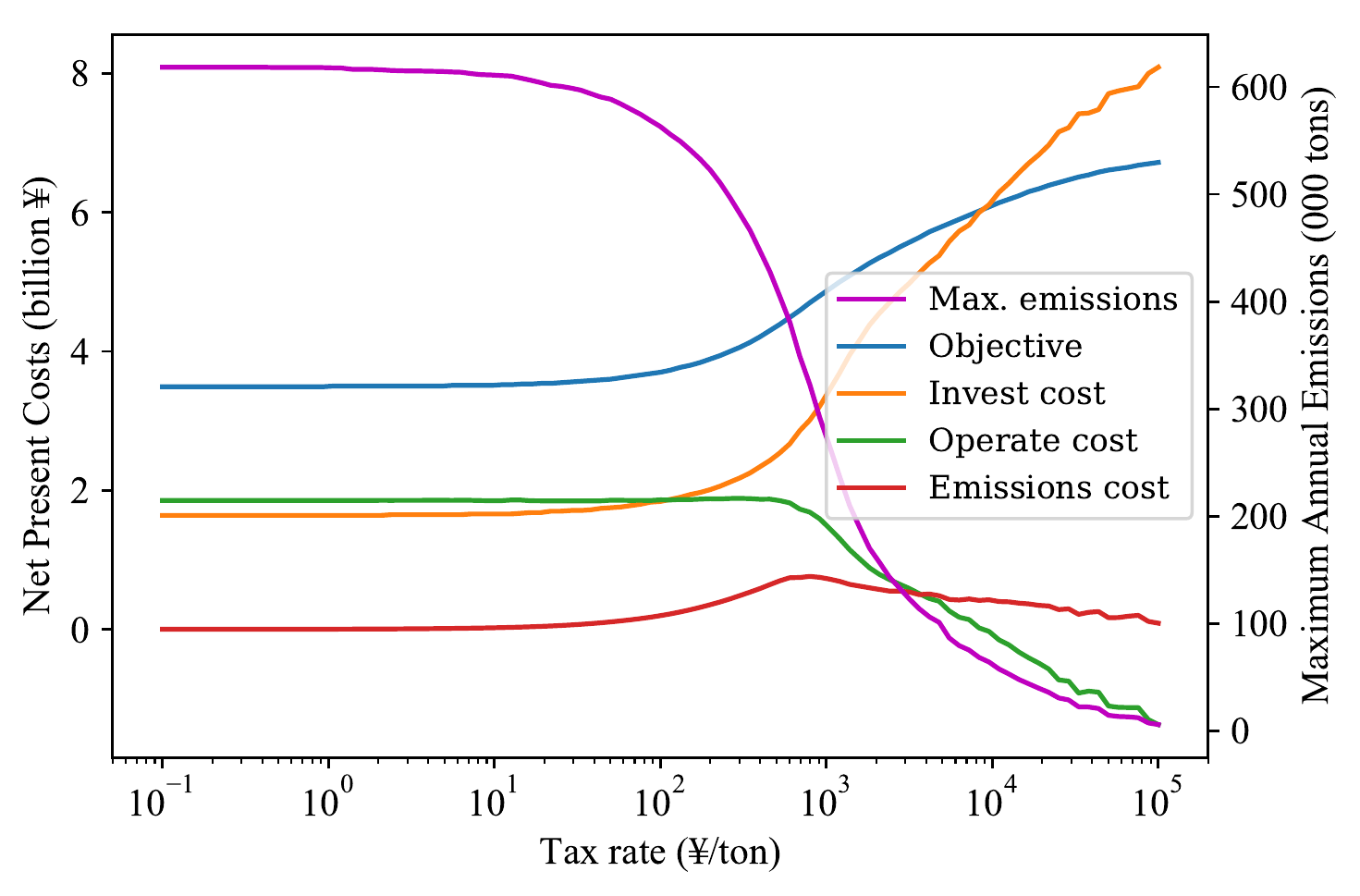}
    \caption{Costs and emissions as a function of carbon tax rate for Framework 2.}
    \label{objective_components_psweep}
\end{figure}

Fig. \ref{objective_components_psweep} shows the cost of each objective component and the resulting maximum emissions as a function of the carbon tax rate for the bi-level regulator/builder-operator framework. The majority of the emissions reductions are achieved at tax rates between 100 and 10,000 \textyen/ton (15-1,500 \$/ton). Table \ref{tax_rate_table} indicates the tax rates required to achieve selected milestones.

\begin{table}[h]
    \centering
    \caption{Framework 2: Tax rates required to meet emissions goals}
    \begin{tabular}{ 
        			>{\centering\arraybackslash}p{0.13\linewidth}
        			>{\centering\arraybackslash}p{0.10\linewidth}
        			>{\centering\arraybackslash}p{0.085\linewidth}
        			>{\centering\arraybackslash}p{0.14\linewidth}
        			>{\centering\arraybackslash}p{0.13\linewidth}
        			>{\centering\arraybackslash}p{0.13\linewidth}
            		}
        \hline \hline
        Emissions reduction target & Tax rate (\textyen/ton) & Total cost (B\textyen) & Investment cost (B\textyen) & Net operate cost (B\textyen) & Emissions cost (B\textyen) \\
        \hline
        None & N/A & 3.49 & 1.64 & 1.85 & 0\\
        5\% & 52 & 3.61 & 1.75 & 1.85 & 0.110 \\
        10\% & 112 & 3.73 & 1.87 & 1.86 & 0.217 \\
        25\% & 358 & 4.15 & 2.28 & 1.88 & 0.548 \\
        50\% & 842 & 4.74 & 3.10 & 1.64 & 0.753 \\
        75\% & 2,109 & 5.35 & 4.59 & 0.763 & 0.548 \\
        90\% & 9,988 & 6.10 & 6.21 & -0.106 & 0.408
    \end{tabular}
    \label{tax_rate_table}
\end{table}

\subsection{Results for Framework 3: Builder / Operator}

Fig. \ref{compare_methods} and Table \ref{bo_vs_sbo_premium} show that when the hub operator cannot be trusted to cooperate to achieve emission reduction targets, the hub must be `overbuilt' as compared to the cooperative scenario of Framework 1. For example, a target reduction of 50\% in Framework 1 is achievable with a total cost of 3.62 B\textyen, of which 2.13 B\textyen\hspace{0.1em} is for equipment purchases. On the other hand, when an operator who cares only about operating costs uses these equipment, the resulting emissions only decrease by 28.5\%. To achieve the 50\% reduction in Framework 3, equipment costing 3.07 B\textyen\hspace{0.1em} are needed, with a resulting total cost of 3.89 B\textyen. In effect, lack of cooperation between the builder and operator are responsible for an overall cost increase of 7\% and an increase of 45\% in construction costs.

\begin{table}[h]
    \centering
    \caption{Comparing costs for equivalent targets in Frameworks 1 and 3}
    \begin{tabular}{ 
        			>{\centering\arraybackslash}p{0.15\linewidth}
        			>{\centering\arraybackslash}p{0.155\linewidth}
        			>{\centering\arraybackslash}p{0.155\linewidth}
        			>{\centering\arraybackslash}p{0.15\linewidth}
        			>{\centering\arraybackslash}p{0.15\linewidth}
            		}
        \hline
        Emissions reduction target & Framework 1 cost (B\textyen) & Framework 3 cost (B\textyen) & Total cost increase & Investment cost increase \\
        \hline \hline
        25\% & 3.52 & 3.57 & 1.3\% & 27.5\% \\
        50\% & 3.62 & 3.89 & 7.3\% & 44.6\% \\
        75\% & 4.35 & 5.30 & 21.8\% & 35.7\% \\
        100\% & 6.80 & 8.14 & 19.6\% & 33.8\%
    \end{tabular}
    \label{bo_vs_sbo_premium}
\end{table}

\subsection{Results for Framework 4: Regulator / Builder / Operator Framework}

When an energy hub builder is forced to consider a social cost of carbon when designing an energy hub, more money is invested in more-expensive but less-polluting equipment. Table \ref{rbo_table} shows the results for several emissions reductions benchmarks. When results from this emissions pricing framework are compared against results from the regulator/builder-operator framework (Table \ref{tax_rate_table}), it can be seen that emissions targets up to 50\% can be realized at a lower total cost. This is due in part to the fact that the builder selects equipment as if emissions were taxed, but the operator doesn't actually have to pay the price for their emissions. At emissions reduction targets of 75\% and higher, however, Framework 4 results in higher costs than Framework 2. Since the operator does not see the price of emissions, the hub must be `overbuilt' to account for the operator's indifference toward emissions. The result is a hub with high PV capacity and no connection to the district heat or gas networks.

\begin{table}[h]
    \centering
    \caption{Tax rates required to meet emissions goals in Framework 4}
    \begin{tabular}{ 
        			>{\centering\arraybackslash}p{0.13\linewidth}
        			>{\centering\arraybackslash}p{0.10\linewidth}
        			>{\centering\arraybackslash}p{0.085\linewidth}
        			>{\centering\arraybackslash}p{0.14\linewidth}
        			>{\centering\arraybackslash}p{0.13\linewidth}
        			>{\centering\arraybackslash}p{0.13\linewidth}
            		}
        \hline \hline
        Emissions reduction target & SCoC rate (\textyen/ton) & Total cost (B\textyen) & Investment cost (B\textyen) & Net operate cost (B\textyen) & Total SCoC (B\textyen) \\
        \hline
        None & N/A & 3.49 & 1.64 & 1.85 & 0\\
        5\% & 5 & 3.50 & 1.79 & 1.71 & 0.034 \\
        10\% & 66 & 3.55 & 1.86 & 1.69 & 0.363 \\
        25\% & 87 & 3.58 & 2.12 & 1.45 & 0.448 \\
        50\% & 346 & 3.92 & 3.20 & 0.722 & 1.242 \\
        75\% & 794 & 5.49 & 6.39 & -0.896 & 0.584 \\
        90\% & 2,754 & 6.31 & 7.76 & -1.44 & 0.572
    \end{tabular}
    \label{rbo_table}
\end{table}

\section{Discussion}

The applicability of a policy framework to reduce carbon emissions from one or more greenfield energy hubs depends on the political-economic structure of the jurisdiction where the energy hubs are to be built. Framework 1 is the cheapest way to reach a given emissions target, but relies on the builder and operator's dedication to meeting the emissions target, during both construction and operation, when alternative investment and operation decisions are cheaper. Framework 2 is more effective at economically incentivizing the construction and operation of a hub able to meet emissions targets, because the cost of emissions is internalized. However, implementing carbon taxes can be politically risky and the total cost of construction and operation can be significantly higher than with Framework 1, especially for mid-range emissions reduction targets. Framework 3 is able to achieve emissions targets without the use of carbon taxes, even when the operator only pursues minimum cost operation, but the overall cost can be significantly higher for aggressive emissions reduction targets. Framework 4 is the most complex, but it enables a mandate that a carbon price be considered in investment decisions without requiring a politically sensitive carbon price to be paid during operation.

Carbon emissions targets are usually pledged on a per-year basis. A constraint on the maximum annual emissions (whether via an explicit constraint or via a strategic tax) may only bind against the emissions in a single year. However, the decisions made in order to meet commitments for the critical year (either infrastructure investment or tax rate) will also reduce emissions in all other years. Optimizing to constrain \textit{lifetime} emissions is an alternate approach that would result in different investment and operation decisions, but is less compatible with the annual emissions pledges which are most common today.

\section{Conclusion}

This paper describes and illustrates four low-carbon energy hub design frameworks, using a new formulation of the network constraints for the greenfield energy hub problem. The single-actor cooperative framework results in the lowest overall cost for any given carbon emissions target, while the bi- and tri-level problems are more expensive due to the lack of cooperation from lower-level actors. Using these frameworks, the impact of various policy decisions on the construction and operation of energy hubs can be investigated. Policy questions which can be aided by these frameworks include:
\begin{itemize}
    \item How should campuses be compensated for electricity exported to the grid, based on their contributions of energy as well as reductions in grid emissions?
    \item What secondary effects can be anticipated in response to investments in new grid-scale electricity and district heat assets, which will change the price and carbon intensity of grid-sourced energy?
    \item To what extent do limitations in regional energy transmission reduce the ability of energy hubs to meet emissions targets, or increase their cost of doing so?
\end{itemize}

Holistically investigating questions such as these can bring society closer to affordable, sustainable power systems.

%END MAIN DOCUMENT
\bibliographystyle{IEEEtran}
\bibliography{bibliography}

\section*{Appendix A. Time-Invariant Matrix Definitions}

\begin{align}
    \boldsymbol{Z}
    :=
    &\begin{bmatrix}
        \boldsymbol{H}_{1} \boldsymbol{A}_{1} & \boldsymbol{H}_{2} \boldsymbol{A}_{2} & \dots & \boldsymbol{H}_{|G|} \boldsymbol{A}_{|G|}
    \end{bmatrix}^{T} \displaybreak[0] \nonumber \\
    A_{p,l} =
    &\begin{cases}
        -1 & \text{if branch } l \text{ flows from port } p \\
        1 & \text{if branch } l \text{ flows to port } p \\
        0 & \text{otherwise}
    \end{cases} \displaybreak[0] \nonumber \\
    H_{p_i,p_j} =
    &\begin{cases}
        1 & \text{if } p_i \text{ is a converter output and } p_i = p_j \\
        \eta & \text{if } p_j \text{ is a converter input corresponding}  \\
        & \text{to output } p_i \text{ with efficiency } \eta \\
        \eta & \text{if } p_j \text{ is a storage output corresponding} \\
        & \text{to input } p_i \\
        \frac{1}{\eta} & \text{if } p_j \text{ is a storage output and } p_i = p_j \\
        1 & \text{if } p_i \text{ is a storage output corresponding} \\
        & \text{to virtual port } p_j \\
        0 & \text{otherwise}
    \end{cases} \displaybreak[0] \nonumber \\
    J_{g,p} = 
    &\begin{cases}
        1 & \text{if port } p \text{ is an input to converter } g \\
        & \text{ or the virtual port to storage } g \\
        0 & \text{otherwise}
    \end{cases} \displaybreak[0] \nonumber \\
    U_{m,l} =
    &\begin{cases}
        1 & \text{if grid source } m \text{ supplies branch } l \\
        0 & \text{otherwise}
    \end{cases} \displaybreak[0] \nonumber \\
    W_{m,l} =
    &\begin{cases}
        1 & \text{if branch } l \text{ supplies end-use } m \\
        0 & \text{otherwise}
    \end{cases} \displaybreak[0] \nonumber \\
    K_{l} =
    &\begin{cases}
        1 & \text{if branch } l \text{ is `real'} \\
        0 & \text{if branch } l \text{ is `virtual'}
    \end{cases} \nonumber
\end{align}

\begin{IEEEbiography}[{\includegraphics[width=1in,height=1.25in,clip,keepaspectratio]{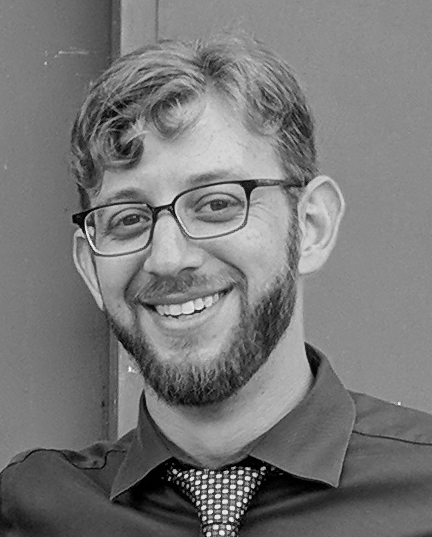}}]{Daniel Olsen} (S'14) received the B.Sc. degree in Mechanical Engineering and Electric Power Engineering from Rensselaer Polytechnic Institute, Troy, NY, USA, in 2010. He is pursuing the Ph.D. degree in Electrical Engineering at the University of Washington, Seattle, WA, USA. Previously, he was a Research Associate with Lawrence Berkeley National Laboratory. His research interests include planning and policies for low-carbon power systems, multiple-energy systems, and distributed flexibility resources.
\end{IEEEbiography}

\begin{IEEEbiography}[{\includegraphics[width=1in,height=1.25in,clip,keepaspectratio]{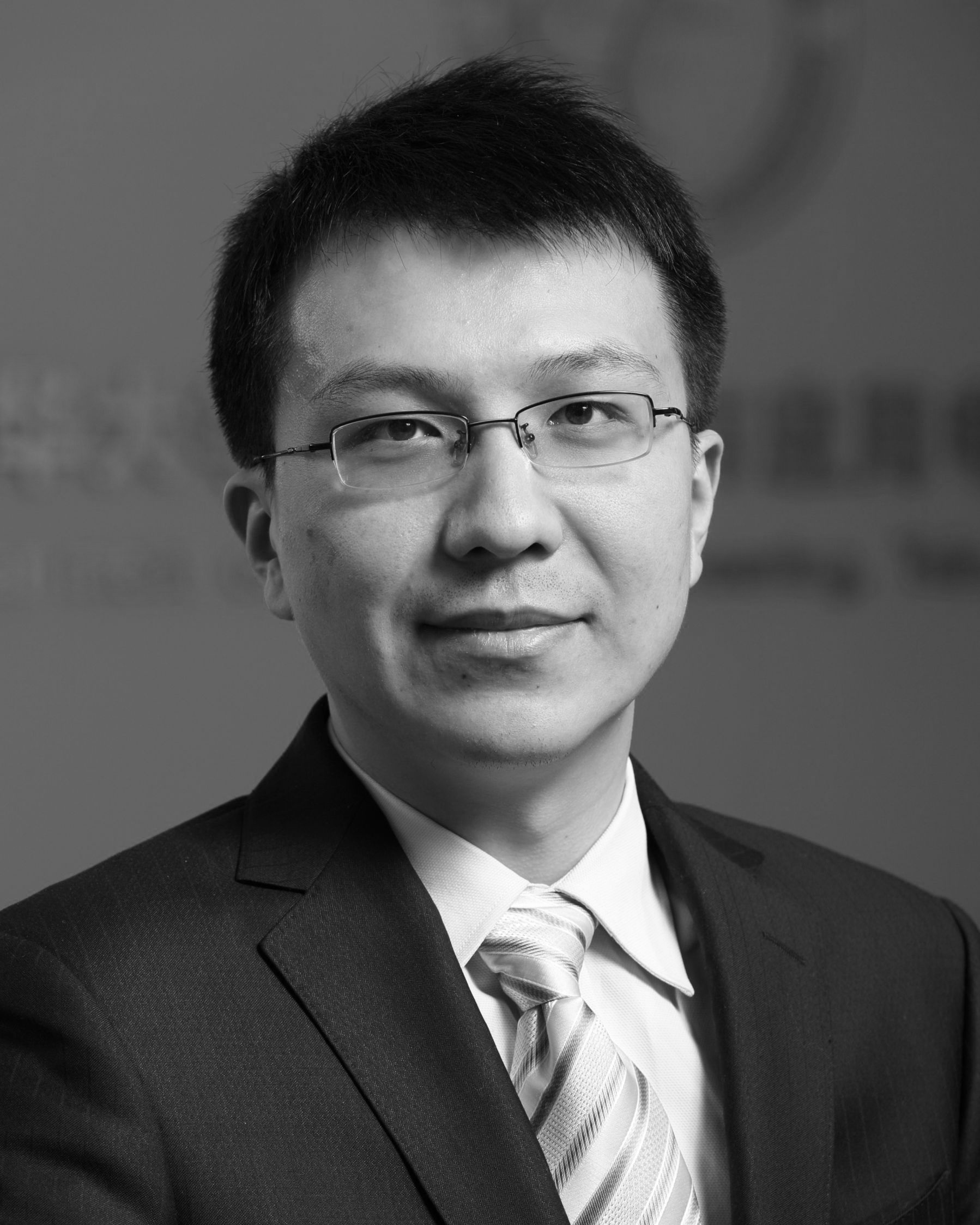}}]{Ning Zhang} (S’10, M’12, SM’18) received both B.S. and Ph.D. from the Electrical Engineering Department of Tsinghua University in China in 2007 and 2012, respectively. 

He is now an Associate Professor at the same university. His research interests include multiple energy system integration, renewable energy, and power system planning and operation.
\end{IEEEbiography}

\begin{IEEEbiography}[{\includegraphics[width=1in,height=1.25in,clip,keepaspectratio]{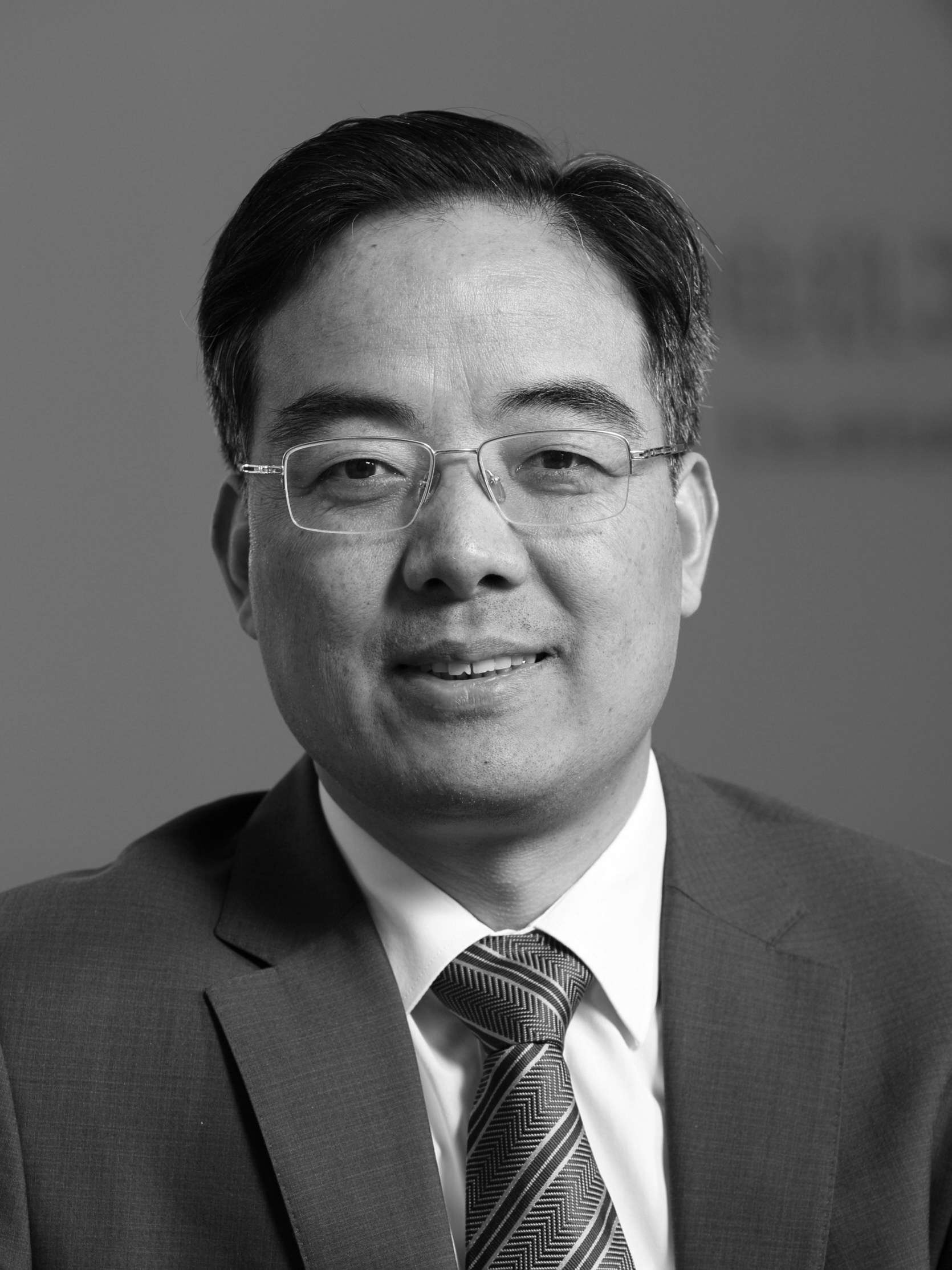}}]{Chongqing Kang} (M’01-SM’07-F’17) received his Ph.D. from the Electrical Engineering Department of Tsinghua University in 1997. 

He is now a Professor at the same university. His research interests include power system planning, power system operation, renewable energy, low carbon electricity technology, and load forecasting.
\end{IEEEbiography}

\begin{IEEEbiography}[{\includegraphics[width=1in,height=1.25in,clip,keepaspectratio]{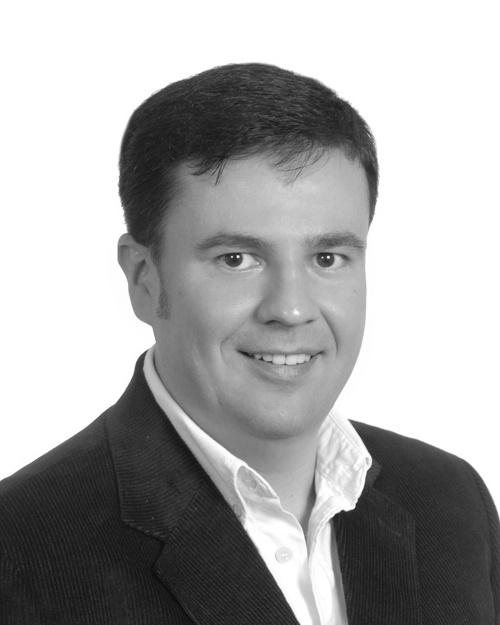}}]{Miguel A. Ortega-Vazquez} (S’97–M’06–SM’14) is a Senior Technical Leader at the Electric Power Research Institute (EPRI) in Palo Alto, CA.  His current research interests include power system operation, power system security, power system economics, integration of renewable energy sources, and the smart grid.  Before joining EPRI, he was Assistant Professor at the University of Washington in Seattle, WA; and before that, Assistant Professor at the Chalmers University of Technology, Sweden.  He holds a Ph.D. from The University of Manchester; M.Sc. from the Universidad Autónoma de Nuevo León, Mexico; and an Electric Engineering degree from the Instituto Tecnológico de Morelia, Mexico.

Dr. Ortega-Vazquez is an Editor of the IEEE Transactions on Smart Grid.
\end{IEEEbiography}

\begin{IEEEbiography}[{\includegraphics[width=1in,height=1.25in,clip,keepaspectratio]{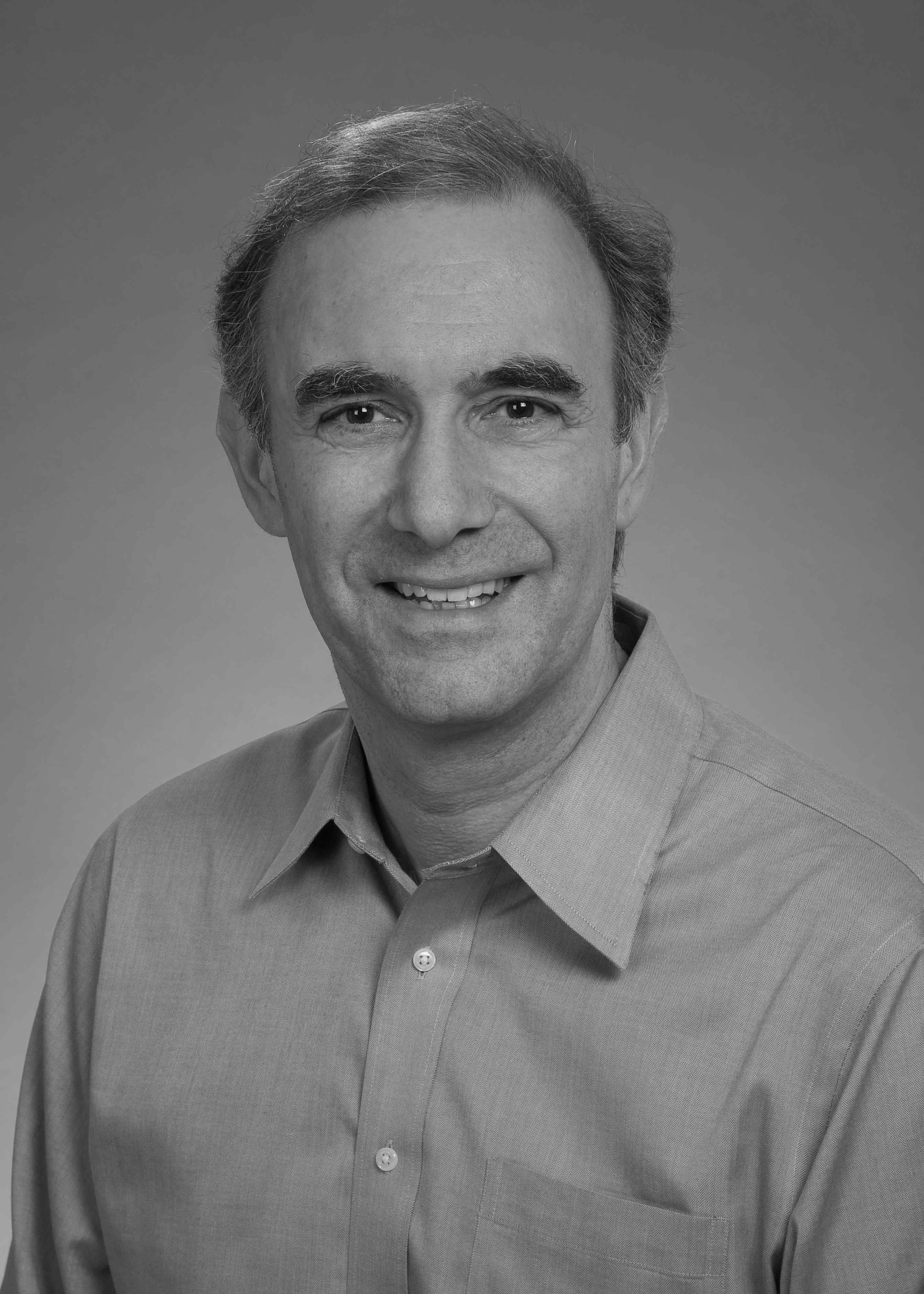}}]{Daniel Kirschen} (M’86–SM’91–F’07) 
is the Donald W. and Ruth Mary Close Professor of Electrical Engineering at the University of Washington. His research focuses on the integration of renewable energy sources in the grid, power system economics and power system resilience. Prior to joining the University of Washington, he taught for 16 years at The University of Manchester (UK). Before becoming an academic, he worked for Control Data and Siemens on the development of application software for utility control centers. He holds a PhD from the University of Wisconsin-Madison and an Electro-Mechanical Engineering degree from the Free University of Brussels (Belgium). He is the author of two books and a Fellow of the IEEE.
\end{IEEEbiography}

\end{document}